\documentclass{amsart}
\usepackage{graphicx}
\usepackage{url}

\theoremstyle{plain}
\newtheorem{theorem}{Theorem}

\newtheorem{conjecture}[theorem]{Conjecture}

\numberwithin{theorem}{section}

\theoremstyle{definition}

\numberwithin{equation}{section}

\newcommand{\R}{{\mathbb R}}
\newcommand{\Z}{{\mathbb Z}}
\newcommand{\C}{{\mathbb C}}

\newcommand{\F}{{\mathbb F}}
\newcommand{\Proj}{{\mathbb P}}
\newcommand{\co}{\colon}
\newcommand{\universal}[2]{U_{#1,\,#2}}

\title[Energy-minimizing point configurations]
{Experimental study of energy-minimizing\\
point configurations on
spheres}

\author[Ballinger]{Brandon Ballinger}
\address{Department of Computer Science\\
University of Washington\\
Seattle, WA 98195} \email{brandonb@cs.washington.edu}

\author[Blekherman]{Grigoriy Blekherman}
\address{Microsoft Research\\
One Microsoft Way\\
Redmond, WA 98052-6399} \email{gblekher@microsoft.com}

\author[Cohn]{Henry Cohn}
\address{Microsoft Research\\
One Microsoft Way\\
Redmond, WA 98052-6399}
\email{cohn@microsoft.com}

\author[Giansiracusa]{Noah Giansiracusa}
\address{Department of Mathematics\\
University of Washington\\
Seattle, WA 98195} \email{noahgian@math.washington.edu}
\curraddr{Department of Mathematics\\
Brown University\\
Providence, RI 02912} \email{noahgian@math.brown.edu}

\author[Kelly]{Elizabeth Kelly}
\address{Department of Mathematics\\
University of Washington\\
Seattle, WA 98195}
\email{beth@math.washington.edu}

\author[Sch\"urmann]{Achill Sch\"urmann}
\address{Institute for Algebra and Geometry\\
University of Magdeburg\\
39106 Magdeburg\\
Germany} \email{achill@math.uni-magdeburg.de}

\thanks{Ballinger and Giansiracusa were supported by
the University of Washington Mathematics Department's NSF VIGRE
grant. Sch\"urmann was supported by the Deutsche
Forschungsgemeinschaft (DFG) under grant SCHU 1503/4-1.}

\date{September 8, 2008}

\begin{document}

\begin{abstract}
In this paper we report on massive computer experiments aimed at
finding spherical point configurations that minimize potential
energy.  We present experimental evidence for two new universal
optima (consisting of $40$ points in $10$ dimensions and $64$
points in $14$ dimensions), as well as evidence that there are no
others with at most $64$ points.  We also describe several other
new polytopes, and we present new geometrical descriptions of some
of the known universal optima.
\end{abstract}

\maketitle

\begin{quote}
[T]he problem of finding the configurations of stable equilibrium
for a number of equal particles acting on each other according to
some law of force\dots is of great interest in connexion with the
relation between the properties of an element and its atomic
weight. Unfortunately the equations which determine the stability
of such a collection of particles increase so rapidly in
complexity with the number of particles that a general
mathematical investigation is scarcely possible. \hfill\break
\phantom{}\ \hfill J.~J.~Thomson, 1897
\end{quote}

\tableofcontents

\section{Introduction} \label{section:intro}

What is the best way to distribute $N$ points over the unit sphere
$S^{n-1}$ in $\R^n$?  Of course the answer depends on the notion
of ``best.''  One particularly interesting case is energy
minimization.  Given a continuous, decreasing function $f \co
(0,4] \to \R$, define the $f$-potential energy of a finite subset
$\mathcal{C} \subset S^{n-1}$ to be
$$
E_f(\mathcal{C}) =
\frac{1}{2}\sum_{\genfrac{}{}{0pt}{}{\scriptstyle x,y \in
\mathcal{C}}{\scriptstyle x \ne y}} f\big(|x-y|^2\big).
$$
(We only need $f$ to be defined on $(0,4]$ because $|x-y|^2 \le 4$
when $|x|^2=|y|^2=1$.  The factor of $1/2$ is chosen for
compatibility with the physics literature, while the use of
squared distance is incompatible but more convenient.) How can one
choose $\mathcal{C} \subset S^{n-1}$ with $|\mathcal{C}|=N$ so as
to minimize $E_f(\mathcal{C})$?  In this paper we report on
lengthy computer searches for configurations with low energy. What
distinguishes our approach from most earlier work on this topic
(see for example \cite{AP-G1, AP-G2, AP-G3, AWRTSDW, AW,
BBCDHNNTW, BCNT1, BCNT2, C, DM, DLT, E, ELSBB, EH, F, GE, HS1,
HS2, H, KaS, KuS, K, LL, M-FMRS, MKS, MDH, P-GDMOD-S, P-GDM, P-GM,
RSZ1, RSZ2, SK, S, T, Wh, Wi}) is that we attempt to treat many
different potential functions on as even a footing as possible.
Much of the mathematical structure of this problem becomes
apparent only when one varies the potential function $f$.
Specifically, we find that many optimal configurations vary in
surprisingly simple, low-dimensional families as $f$ varies.

The most striking case is when such a family is a single point: in
other words, when the optimum is independent of $f$.  Cohn and
Kumar \cite{CK1} defined a configuration to be \textit{universally
optimal\/} if it minimizes $E_f$ for all completely monotonic $f$
(i.e., $f$ is infinitely differentiable and $(-1)^k f^{(k)}(x) \ge
0$ for all $k \ge 0$ and  $x \in (0,4)$, as is the case for
inverse power laws). They were able to prove universal optimality
only for certain very special arrangements.  One of our primary
goals in this paper is to investigate how common universal
optimality is. Was the limited list of examples in \cite{CK1} an
artifact of the proof techniques or a sign that these
configurations are genuinely rare?

Every universally optimal configuration is an optimal spherical
code, in the sense that it maximizes the minimal distance between
the points.  (Consider an inverse power law $f(r) = 1/r^s$.  If
there were a configuration with a larger minimal distance, then
its $f$-potential energy would be lower when $s$ is sufficiently
large.)  However, universal optimality is a far stronger condition
than optimality as a spherical code.  There are optimal spherical
codes of each size in each dimension, but they are rarely
universally optimal. In three dimensions, the only examples are a
single point, two antipodal points, an equilateral triangle on the
equator, or the vertices of a regular tetrahedron, octahedron, or
icosahedron. Universal optimality was proved in \cite{CK1},
building on previous work by Yudin, Kolushov, and Andreev
\cite{Y,KY1,KY2,A1,A2}, and the completeness of this list follows
from a classification theorem due to Leech \cite{L}. See
\cite{CK1} for more details.

In higher dimensions much less is known.  Cohn and Kumar's main
theorem provides a general criterion from which they deduced the
universal optimality of a number of previously studied
configurations.  Specifically, they proved that every spherical
$(2m-1)$-design in which only $m$ distances occur between distinct
points is universally optimal.  Recall that a spherical $d$-design
in $S^{n-1}$ is a finite subset $\mathcal{C}$ of $S^{n-1}$ such
that every polynomial on $\R^n$ of total degree $d$ has the same
average over $\mathcal{C}$ as over the entire sphere.  This
criterion holds for every known universal optimum except one case,
namely the regular $600$-cell in $\R^4$ (i.e., the $H_4$ root
system), for which Cohn and Kumar proved universal optimality by a
special argument.

\begin{table}
\begin{center}
\begin{tabular}{cccc}
$n$ & $N$ & $t$ & Description\\
\hline $2$ & $N$ & $\cos (2\pi /N)$ & $N$-gon\\
$n$ & $N \le n+1$ &  $-1/(N-1)$ & simplex\\
$n$ & $2n$ &  $0$ & cross polytope\\
$3$ & $12$ &  $1/\sqrt{5}$ & icosahedron\\
$4$ & $120$ &  $(1 + \sqrt{5})/4$ & regular $600$-cell\\
$5$ & $16$ &  $1/5$ & hemicube/Clebsch graph\\
$6$ & $27$ &  $1/4$ & Schl\"afli graph/isotropic subspaces\\
$7$ & $56$ &  $1/3$ & equiangular lines\\
$8$ & $240$ &  $1/2$ & $E_8$ root system\\
$21$ & $112$ & $1/9$ & isotropic subspaces\\
$21$ & $162$  & $1/7$ & $(162, 56, 10, 24)$ strongly regular graph\\
$22$ & $100$  & $1/11$ & Higman-Sims graph\\
$22$ & $275$  & $1/6$ & McLaughlin graph\\
$22$ & $891$  & $1/4$ & isotropic subspaces\\
$23$ & $552$  & $1/5$ & equiangular lines\\
$23$ & $4600$  & $1/3$ & kissing configuration of the following\\
$24$ & $196560$ &  $1/2$ & Leech lattice minimal vectors\\
$q\frac{q^3+1}{q+1}$ & $(q+1)(q^3+1)$ & $ 1/q^2$ &
isotropic subspaces ($q$ is a prime power)\\
\\ % Bad to skip line?
\end{tabular}
\end{center}
\caption{The known universal optima.}\label{table:universal}
\end{table}

A list of all known universal optima is given in
Table~\ref{table:universal}.  Here $n$ is the dimension of the
Euclidean space, $N$ is the number of points, and $t$ is the
greatest inner product between distinct points in the
configuration (i.e., the cosine of the minimal angle).  For
detailed descriptions of these configurations, see Section~1 of
\cite{CK1}.  Each is uniquely determined by the parameters listed
in Table~\ref{table:universal}, except for the configurations
listed on the last line.  For that case, when $q=p^\ell$ with $p$
an odd prime, there are at least $\lfloor (\ell-1)/2 \rfloor$
distinct universal optima (see \cite{CGS} and \cite{Ka}).
Classifying these optima is equivalent to classifying generalized
quadrangles with parameters $(q,q^2)$, which is a difficult
problem in combinatorics.  In the other cases from
Table~\ref{table:universal}, when uniqueness holds, we use the
notation $\universal{N}{n}$ for the unique $N$-point universal
optimum in $\R^n$.

Each of the configurations in Table~\ref{table:universal} had been
studied before it appeared in \cite{CK1}, and was already known to
be an optimal spherical code. In fact, when $N \ge 2n+1$ and $n >
4$, the codes on this list are exactly those that have been proved
optimal. Cohn and Kumar were unable to resolve the question of
whether Table~\ref{table:universal} is the complete list of
universally optimal codes, except when $n \le 3$.  All that is
known in general is that any new universal optimum must have $N
\ge 2n+1$ (Proposition~1.4 in \cite{CK1}). It does not seem
plausible that the current list is complete, but it is far from
obvious where to find any others.

Each known universal optimum is a beautiful mathematical object,
connected to various important exceptional structures (such as
special lattices or groups).  Our long-term hope is to develop
automated tools that will help uncover more such objects.  In this
paper we do not discover any configurations as fundamental as
those in Table~\ref{table:universal}, but perhaps our work is a
first step in that direction.

Table~\ref{table:universal} shows several noteworthy features.
When $n \le 4$, the codes listed are the vertices of regular
polytopes, specifically those with simplicial facets. When $5 \le
n \le 8$, the table also includes certain semiregular polytopes
(their facets are simplices and cross polytopes, with two cross
polytopes and one simplex arranged around each $(n-3)$-dimensional
face). The corresponding spherical codes are all affine cross
sections of the minimal vectors of the $E_8$ root lattice.
Remarkably, no universal optima are known for $9 \le n \le 20$,
except for the simplices and cross polytopes, which exist in every
dimension. This gap is troubling---why should these dimensions be
disfavored? For $21 \le n \le 24$ nontrivial universal optima are
known; they are all affine cross sections of the minimal vectors
of the Leech lattice (and are no longer the vertices of
semiregular polytopes). Finally, in high dimensions a single
infinite sequence of nontrivial universal optima is known.

It is not clear how to interpret this list.  For example, is the
dimension gap real, or merely an artifact of humanity's limited
imagination?  One of our conclusions in this paper is that
Table~\ref{table:universal} is very likely incomplete but appears
closer to complete than one might expect.

\subsection{Experimental results}

\begin{table}
\begin{center}
\begin{tabular}{cccc}
$n$ & $N$ & $t$ & References\\
\hline
$10$ & $40$ & $1/6$ & Conway, Sloane, and Smith \cite{S}, Hovinga \cite{H}\\
$14$ & $64$ & $1/7$ & Nordstrom and Robinson \cite{NR},
de Caen and van Dam \cite{dCvD},\\
& & & Ericson and Zinoviev \cite{EZ}\\
\\ % Bad to skip line?
\end{tabular}
\end{center}
\caption{New conjectured universal
optima.}\label{table:conjectures}
\end{table}

One outcome of our computer searches is two candidate universal
optima, listed in Table~\ref{table:conjectures} and described in
more detail in Section~\ref{section:newunivopt}. These
configurations were located through massive computer searches: for
each of many pairs $(n,N)$, we repeatedly picked $N$ random points
on $S^{n-1}$ and performed gradient descent to minimize potential
energy.  We focused on the potential function $f(r) =
1/r^{n/2-1}$, because $x \mapsto 1/|x|^{n-2}$ is the unique
nonconstant radial harmonic function on $\R^n\setminus \{0\}$, up
to scalar multiplication (recall that distance is squared in the
definition of $E_f$). When $n=3$, the $f$-potential energy for
this function $f$ is the Coulomb potential energy from
electrostatics, and this special case has been extensively studied
by mathematicians and other scientists. In higher dimensions, this
potential function has frequently been studied as a natural
generalization of electrostatics; we call it the harmonic
potential function.

\begin{figure}
\begin{center}
\includegraphics{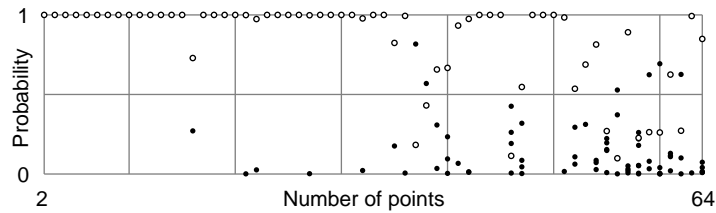}
\end{center}
\caption{Probabilities of local minima for harmonic energy in
$\R^3$ (based on $1000$ trials).  White circles denote the
conjectured harmonic optima.} \label{fig:probs}
\end{figure}

Because there are typically numerous local minima for harmonic
energy, we repeated this optimization procedure many times with
the hope of finding the global minimum.  For low numbers of points
in low dimensions, the apparent global minimum occurs fairly
frequently. Figure~\ref{fig:probs} shows data from three
dimensions.  In higher dimensions, there are usually more local
minima and the true optimum can occur very infrequently.

For each conjectured optimum for harmonic energy, we attempted to
determine whether it could be universally optimal. We first
determined whether it is in equilibrium under all possible force
laws (i.e., ``balanced'' in the terminology of Leech \cite{L}).
That holds if and only if for each point $x$ in the configuration
and each distance $d$, the sum of all points in the code at
distance $d$ from $x$ is a scalar multiple of $x$.  If this
criterion fails, then there is some inverse power law under which
the code is not even in equilibrium, let alone globally minimal,
so it cannot possibly be universally optimal. Most of the time,
the code with the lowest harmonic potential energy is not
balanced. When it is balanced, we compared several potential
functions to see whether we could disprove universal optimality.
By Theorem~9b in \cite[p.~154]{Wid}, it suffices to look at the
potential functions $f(r) = (4-r)^k$ with $k \in \{0,1,2,\dots\}$
(on each compact subinterval of $(0,4]$, every completely
monotonic function can be approximated arbitrarily closely by
positive linear combinations of these potential functions).
Because these functions do not blow up at $r=0$, numerical
calculations with them often converge more slowly than they do for
inverse power laws (nearby points can experience only a bounded
force pushing them apart), so they are not a convenient choice for
initial experimentation. However, they play a fundamental role in
detecting universal optima.

%update me
To date, our search has led us to $58$ balanced configurations
with at most $64$ points (and at least $2n+1$ in dimension $n$)
that appear to minimize harmonic energy and were not already known
to be universally optimal.  In all but two cases, we were able to
disprove universal optimality, but the remaining two cases (those
listed in Table~\ref{table:conjectures}) are unresolved.  We
conjecture that they are in fact universally optimal.

\begin{figure}
\begin{center}
\includegraphics{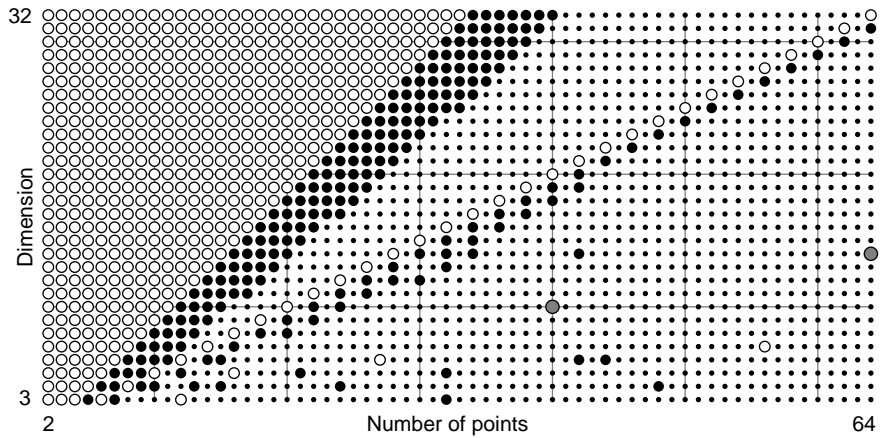}
\end{center}
\caption{Status of conjectured harmonic optima with up to $64$
points in at most $32$ dimensions: white circle denotes universal
optimum, large gray circle denotes conjectured universal optimum,
black circle denotes balanced configuration that is not
universally optimal, tiny black circle denotes unbalanced
configuration.} \label{fig:balanced}
\end{figure}

Figure~\ref{fig:balanced} presents a graphical overview of our
data.  The triangle of white circles on the upper left represents
the simplices, and the diagonal line of white circles represents
the cross polytopes. Between them, one can see that the pattern is
fairly regular, but as one moves right from the cross polytopes
all structure rapidly vanishes.  There is little hope of finding a
simple method to predict where balanced harmonic optima can be
found, let alone universal optima.  It also does not seem likely
that general universal optima can be characterized by any variant
of Cohn and Kumar's criterion.

Besides the isotropic subspace universal optima from
Table~\ref{table:universal} and the other universal optima with
the same parameters, we can conjecture only one infinite family of
balanced harmonic optima with more than $2n$ points in $\R^n$,
namely the diplo-simplices with $2n+2$ points in $\R^n$ for $n \ge
6$ (see Subsection~\ref{subsec:diplo}).  Certainly no others are
apparent in Figure~\ref{fig:balanced}, but the isotropic subspace
optima from Table~\ref{table:universal} are sufficiently large and
exotic that it would be foolish to conjecture that there are no
other infinite families.

\begin{table}
\begin{center}
\begin{tabular}{ccc}
$n$ & $N$ & $t$\\
\hline
$3$ & $32$ & $\sqrt{75+30\sqrt{5}}/15$\\
$4$ & $10$ &  $1/6$ or $(\sqrt{5}-1)/4$\\
% Note negative space:
$4$ & $13$ & $\big(\!\cos(4\pi/13) + \cos(6\pi/13)\big)/2$\\
$4$ & $15$ &  $1/\sqrt{8}$\\
$4$ & $24$ &  $1/2$\\
$4$ & $48$ &  $1/\sqrt{2}$\\
$5$ & $21$ &  $1/\sqrt{10}$\\
$5$ & $32$ &  $1/\sqrt{5}$\\
$n \ge 6$ & $2n+2$ & $1/n$\\
$6$ & $42$ & $2/5$\\
$6$ & $44$ & $1/\sqrt{6}$\\
$6$ & $126$ & $\sqrt{3/8}$\\
$7$ & $78$ & $3/7$\\
$7$ & $148$ & $\sqrt{2/7}$\\
$8$ & $72$ & $5/14$\\
$9$ & $96$ & $1/3$\\
$14$ & $42$ & $1/10$\\
$16$ & $256$ & $1/4$\\
\\ % Bad to skip line?
\end{tabular}
\end{center}
\caption{Conjectured harmonic optima that are balanced,
irreducible, and not universally optimal (see
Section~\ref{section:balanced} for descriptions).}
\label{table:balanced}
\end{table}

\begin{table}
\begin{center}
\begin{tabular}{ccc}
$n$ & $N$ & $t$\\
\hline
$7$ & $182$ & $1/\sqrt{3}$\\
$15$ & $128$ & $1/5$\\
\\ % Bad to skip line?
\end{tabular}
\end{center}
\caption{Unresolved conjectured harmonic optima.}
\label{table:unresolved}
\end{table}

Table~\ref{table:balanced} lists the cases in which we found a
balanced harmonic optimum but were able to disprove universal
optimality, with one systematic exception: we omit configurations
that are reducible, in the sense of being unions of orthogonal,
lower-dimensional configurations (this terminology is borrowed
from the theory of root systems). Reducible configurations are in
no sense less interesting or fruitful than irreducible ones.  For
example, cross polytopes can be reduced all the way to
one-dimensional pieces. However, including reducible
configurations would substantially lengthen
Table~\ref{table:balanced} without adding much more geometrical
content.

Table~\ref{table:unresolved} lists two more unresolved cases. They
both appear to be harmonic optima and are balanced, and we have
not been able to prove or disprove universal optimality.  Unlike
the two cases in Table~\ref{table:conjectures}, we do not
conjecture that they are universally optimal, because each is
closely analogous to a case in which universal optimality fails
($182$ points in $\R^7$ is analogous to $126$ points in $\R^6$,
and $128$ points in $\R^{15}$ is analogous to $256$ points in
$\R^{16}$).  On the other hand, each is also analogous to a
configuration we know or believe is universally optimal ($240$
points in $\R^8$ and $64$ points in $\R^{14}$, respectively). We
have not been able to disprove universal optimality in the cases
in Table~\ref{table:unresolved}, but they are sufficiently large
that our failure provides little evidence in favor of universal
optimality.

Note that the data presented in Tables~\ref{table:balanced}
and~\ref{table:unresolved} may not specify the configurations
uniquely.  For example, for $48$ points in $\R^4$ there is a
positive-dimensional family of configurations with maximal inner
product $1/\sqrt{2}$ (which is not the best possible value,
according to Sloane's tables \cite{S}). See
Section~\ref{section:balanced} for explicit constructions of the
conjectured harmonic optima.

It is worth observing that several famous configurations do not
appear in Tables~\ref{table:balanced} or~\ref{table:unresolved}.
Most notably, the cubes in $\R^n$ with $n \ge 3$, the
dodecahedron, the $120$-cell, and the $D_5$, $E_6$, and $E_7$ root
systems are suboptimal for harmonic energy.  Many of these
configurations have more than $64$ points, but we have included in
the tables all configurations we have analyzed, regardless of
size.

In each case listed in Tables~\ref{table:balanced}
and~\ref{table:unresolved}, our computer programs returned
floating point approximations to the coordinates of the points in
the code, but we have been able to recognize the underlying
structure exactly.  That is possible largely because these codes
are highly symmetric, and once one has uncovered the symmetries
the remaining structure is greatly constrained. By contrast, for
most numbers of points in most dimensions, we cannot even
recognize the minimal harmonic energy as an exact algebraic number
(although it must be algebraic, because it is definable in the
first-order theory of the real numbers).

\subsection{New universal optima}

Both codes listed in Table~\ref{table:conjectures} have been
studied before. The first code was discovered by Conway, Sloane,
and Smith \cite{S} as a conjecture for an optimal spherical code
(and discovered independently by Hovinga \cite{H}). The second can
be derived from the Nordstrom-Robinson binary code \cite{NR} or as
a spectral embedding of an association scheme discovered by de
Caen and van Dam \cite{dCvD} (take $t=1$ in Theorem~2 and
Proposition~7(i) in \cite{dCvD} and then project the standard
orthonormal basis into a common eigenspace of the operators in the
Bose-Mesner algebra of the association scheme). We describe both
codes in greater detail in Section~\ref{section:newunivopt}.

Neither code satisfies the condition from \cite{CK1} for universal
optimality: both are spherical $3$-designs (but not $4$-designs),
with four distances between distinct points in the $40$-point code
and three in the $64$-point code.  That leaves open the
possibility of an ad hoc proof, similar to the one Cohn and Kumar
gave for the regular $600$-cell, but the techniques from
\cite{CK1} do not apply.

To test universal optimality, we have carried out $1000$ random
trials with the potential function $f(r) = (4-r)^k$ for each $k$
from $1$ to $25$. We have also carried out $1000$ trials using
Hardin and Sloane's program Gosset \cite{HSl} to construct good
spherical codes (to take care of the case when $k$ is large). Of
course these experimental tests fall far short of a rigorous
proof, but the codes certainly appear to be universally optimal.

We believe that they are the only possible new universal optima
consisting of at most $64$ points, because we have searched the
space of such codes fairly thoroughly.  By Proposition~1.4 in
\cite{CK1}, any new universal optimum in $\R^n$ must contain at
least $2n+1$ points.  There are $812$ such cases with at most $64$
points in dimension at least $4$.  In each case, we have completed
at least $1000$ random trials (and usually more).  There is no
guarantee that we have found the global optimum in any of these
cases, because it could have a tiny basin of attraction.  However,
a simple calculation shows that it is 99.99\% likely that in every
case we have found every local minimum that occurs at least 2\% of
the time.  We have probably not always found the true optimum, but
we believe that we have found every universal optimum within the
range we have searched.

% UPDATE

We have made our tables of conjectured harmonic optima for up to
$64$ points in up to $32$ dimensions available via the world wide
web at
\begin{center}
\url{http://aimath.org/data/paper/BBCGKS2006/}.
% UPDATE ME
\end{center}
They list the best energies we have found and the coordinates of
the configurations that achieve them.  We would be grateful for
any improvements, and we intend to keep the tables up to date,
with credit for any contributions received from others.

In addition to carrying out our own searches for universal optima,
we have examined Sloane's tables \cite{S} of the best spherical
codes known with at most $130$ points in $\R^4$ and $\R^5$, and we
have verified that they contain no new universal optima.  We
strongly suspect that there are no undiscovered universal optima
of any size in $\R^4$ or $\R^5$, based on Sloane's calculations as
well as our searches, but it would be difficult to give definitive
experimental evidence for such an assertion (we see no convincing
arguments for why huge universal optima should not exist).

In general, our searches among larger codes have been far less
exhaustive than those up to $64$ points: we have at least briefly
examined well over four thousand different pairs $(n,N)$, but
generally not in sufficient depth to make a compelling case that
we have found the global minimum.  (Every time we found a balanced
harmonic optimum, with the exception of $128$ points in $\R^{15}$
and $256$ points in $\R^{16}$,
% update me
we completed at least $1000$ trials to test whether it was really
optimal. However, we have not completed nearly as many trials in
most other cases, and in any case $1000$ trials is not enough when
studying large configurations.) Nevertheless, our strong
impression is that universal optima are rare, and certainly that
there are few small universal optima with large basins of
attraction.

\section{Methodology}

\subsection{Techniques}

As discussed in the introduction, to minimize potential energy we
apply gradient descent, starting from many random initial
configurations. That is an unsophisticated approach, because
gradient descent is known to perform more slowly in many
situations than competing methods such as the conjugate gradient
algorithm. However, it has performed adequately in our
computations. Furthermore, gradient descent has particularly
intuitive dynamics. Imagine particles immersed in a medium with
enough viscosity that they never build up momentum.  When a force
acts on them according to the potential function, the
configuration undergoes gradient descent.  By contrast, for most
other optimization methods the motion of the particles is more
obscure, so for example it is more difficult to interpret
information such as sizes of basins of attraction.

Once one has approximate coordinates, one can use the multivariate
analogue of Newton's method to compute them to high precision (by
searching for a zero of the gradient vector).  Usually we do not
need to do this, because the results of gradient descent are
accurate enough for our purposes, but it is a useful tool to have
available.

\begin{figure}
\begin{center}
\includegraphics{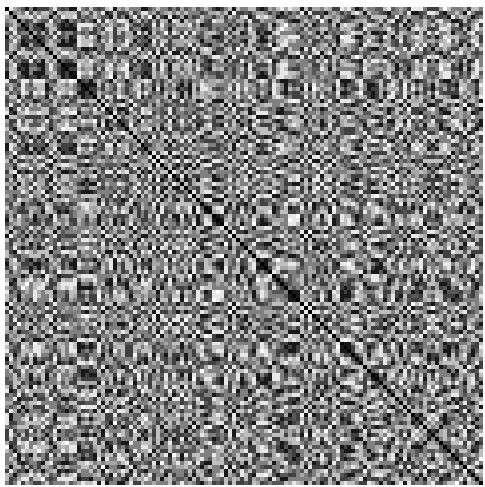}
\end{center}
\caption{The Gram matrix for a regular $600$-cell (black denotes
$1$, white denotes $-1$, and gray interpolates between them), with
the points ordered as returned by our gradient descent software.}
\label{fig:120gram}
\end{figure}

\begin{figure}
\begin{center}
\includegraphics{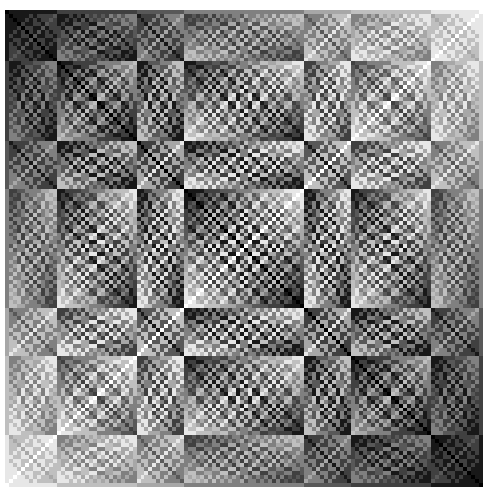}
\end{center}
\caption{The Gram matrix for a regular $600$-cell, with the points
ordered so as to display structure.} \label{fig:120gram2}
\end{figure}

Obtaining coordinates is simply the beginning of our analysis.
Because the coordinates encode not only the relative positions of
the points but also an arbitrary orthogonal transformation of the
configuration, interpreting the data can be subtle.  A first step
is to compute the Gram matrix.  In other words, given points
$x_1,\dots,x_N$, compute the $N \times N$ matrix $G$ whose entries
are given by $G_{i,j} = \langle x_i, x_j \rangle$. The Gram matrix
is invariant under orthogonal transformations, so it encodes
almost precisely the information we care about.  Its only drawback
is that it depends on the arbitrary choice of how the points are
ordered.  That may sound like a mild problem, but there are many
permutations of the points and it is far from clear how to choose
one that best exhibits the configuration's underlying structure:
compare Figure~\ref{fig:120gram} with Figure~\ref{fig:120gram2}.

With luck, one can recognize the entries of the Gram matrix as
exact algebraic numbers: more frequently than one might expect,
they are rational or quadratic irrationals.  Once one specifies
the entire Gram matrix, the configuration is completely
determined, up to orthogonal transformations.  Furthermore, one
can easily prove that the configuration exists (keep in mind that
it may not be obvious that there actually is such an arrangement
of points, because it was arrived at via inexact calculations). To
do so, one need only check that the Gram matrix is symmetric, it
is positive semidefinite, and its rank is at most $n$.  Every such
matrix is the Gram matrix of a set of $N$ points in $\R^n$, and if
the diagonal entries are all $1$ then the points lie on $S^{n-1}$.

\begin{figure}
\begin{center}
\includegraphics{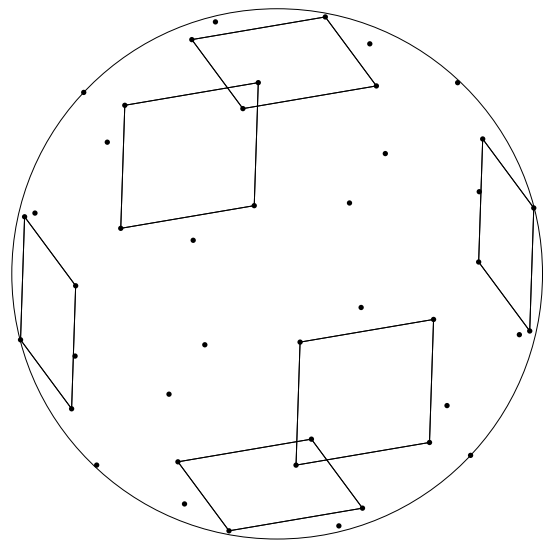}
\end{center}
\caption{An orthogonal projection of the conjectured harmonic
optimum with $44$ points in $\R^3$ onto a random plane.  Line
segments connect points at the minimal distance.} \label{fig:3-44}
\end{figure}

\begin{figure}
\begin{center}
\includegraphics{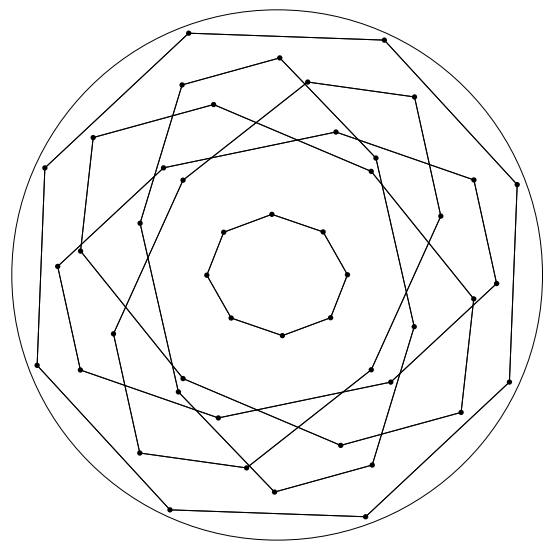}
\end{center}
\caption{An orthogonal projection of the conjectured harmonic
optimum with $48$ points in $\R^4$ onto a random plane.  Line
segments connect points at the minimal distance.} \label{fig:4-48}
\end{figure}

Unfortunately, the exact Gram matrix entries are not always
apparent from the numerical data.  There is also a deeper reason
why simply recognizing the Gram matrix is unsatisfying: it
provides only the most ``bare bones'' description of the
configuration.  Many properties, such as symmetry or connections
to other mathematical structures, are far from apparent given only
the Gram matrix, as one can see from Figures~\ref{fig:120gram}
and~\ref{fig:120gram2}.

Choosing the right method to visualize the data can make the
underlying patterns clearer.  For example, projections onto
low-dimensional subspaces are often illuminating.  Determining the
most revealing projection can itself be difficult, but sometimes
even a random projection sheds light on the structure.  For
example, Figures~\ref{fig:3-44} and~\ref{fig:4-48} are projections
of the harmonic optima with $44$ points in $\R^3$ and $48$ points
in $\R^4$, respectively, onto random planes.  The circular outline
is the boundary of the projection of the sphere, and the line
segments pair up points separated by the minimal distance.
Figure~\ref{fig:3-44} shows a disassembled cube (in a manner
described later in this section), while Figure~\ref{fig:4-48} is
made up of octagons (see Subsection~\ref{subsec:48in4} for a
description).

The next step in the analysis is the computation of the
automorphism group.  In general that is a difficult task, but we
can make use of the wonderful software Nauty written by McKay
\cite{M}. Nauty can compute the automorphism group of a graph as a
permutation group on the set of vertices;  more generally, it can
compute the automorphism group of a vertex-labeled graph.  We make
use of it as follows.

Define a combinatorial automorphism of a configuration to be a
permutation of the points that preserves inner products
(equivalently, distances).  If one forms an edge-labeled graph by
placing an edge between each pair of points, labeled by their
inner product, then the combinatorial automorphism group is the
automorphism group of this labeled graph.  Nauty is not directly
capable of computing such a group, but it is straightforward to
reduce the problem to that of computing the automorphism group of
a related vertex-labeled graph.  Thus, one can use Nauty to
compute the combinatorial automorphism group.

Fortunately, combinatorial automorphisms are the same as geometric
symmetries, provided the configuration spans $\R^n$. Specifically,
every combinatorial automorphism is induced by a unique orthogonal
transformation of $\R^n$.  (When the points do not span $\R^n$,
the orthogonal transformations are not unique, because there are
nontrivial orthogonal transformations that fix the subspace
spanned by the configuration.) Thus, Nauty provides an efficient
method for computing the symmetry group.

Unfortunately, it is difficult to be certain that one has computed
the correct group.  Two inner products that appear equal
numerically may differ by a tiny amount, in which case the
computed symmetry group may be too large.  However, that is rarely
a problem even with single-precision floating point arithmetic,
and it is difficult to imagine a fake symmetry that appears real
to one hundred decimal places.

Once the symmetry group has been obtained, many further questions
naturally present themselves.  Can one recognize the symmetry
group as a familiar group?  How does its representation on $\R^n$
break up into irreducibles?  What are the orbits of its action on
the configuration?

Analyzing the symmetries of the configuration frequently
determines much of the structure, but usually not all of it.  For
example, consider the simplest nontrivial case, namely five points
on $S^2$.  There are two natural ways to arrange them: with two
antipodal points and three points forming an equilateral triangle
on the orthogonal plane between them, or as a pyramid with four
points forming a square in the hemisphere opposite a single point
(and equidistant from it).  In the first case everything is
determined by the symmetries, but in the second there is one free
parameter, namely how far the square is from the point opposite
it. As one varies the potential function, the energy-minimizing
value of this parameter will vary.  (We conjecture that for every
completely monotonic potential function, one of the configurations
described in this paragraph globally minimizes the energy, but we
cannot prove it.)

We define the \textit{parameter count\/} of a configuration to be
the dimension of the space of nearby configurations that can be
obtained from it by applying an arbitrary radial force law between
all pairs of particles. For example, balanced configurations are
those with zero parameters, and the family with a square opposite
a point has one parameter.

To compute the parameter count for an $N$-point configuration,
start by viewing it as an element of $(S^{n-1})^N$ (by ordering
the points).  Within the tangent space of this manifold, for each
radial force law there is a tangent vector. To form a basis for
all these force vectors, look at all distances $d$ that occur in
the configuration, and for each of them consider the tangent
vector that pushes each pair of points at distance $d$ in opposite
directions but has no other effects. All force vectors are linear
combinations of these ones, and the dimension of the space they
span is the parameter count for the configuration. (One must be
careful to use sufficiently high-precision arithmetic, as when
computing the symmetry group.)

This information is useful because in a sense it shows how much
humanly understandable structure we can expect to find.  For
example, in the five-point configuration with a square opposite a
point, the distance between them will typically be some
complicated number depending on the potential function.  In
principle one can describe it exactly, but in practice it is most
pleasant to treat it as a black box and describe all the other
distances in the configuration in terms of it.  The parameter
count tells how many independent parameters one should expect to
arrive at.  When the count is zero or one, it is reasonable to
search for an elegant description, whereas when the count is
twenty, it is likely that the configuration is unavoidably
complex.

\begin{figure}
\begin{center}
\includegraphics{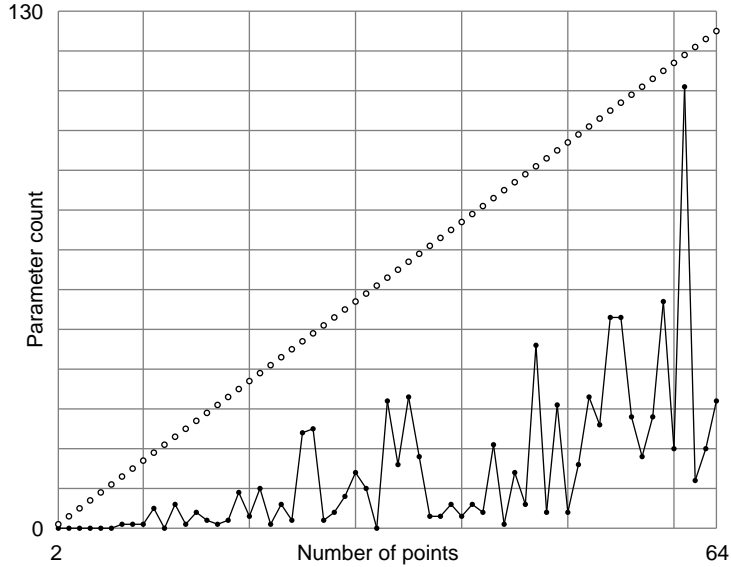}
\end{center}
\caption{Parameter counts for conjectured harmonic optima in
$\R^3$. Horizontal or vertical lines occur at multiples of ten,
and white circles denote the dimension of the configuration
space.} \label{fig:params}
\end{figure}

Figure~\ref{fig:params} shows the parameter counts of the
conjectured harmonic optima in $\R^3$ with at most $64$ points,
compared with the dimension of the full space of all
configurations of their size. The counts vary wildly but are often
quite a bit smaller than one might expect.  Two striking examples
are $61$ points with $111$ parameters, for which there is likely
no humanly understandable description, and $44$ points with one
parameter. The $44$-point configuration consists of the vertices
of a cube and centers of its edges together with the $24$-point
orbit (under the cube's symmetry group) of a point on a diagonal
of a face, all projected onto a common sphere.  The optimal choice
of the point on the diagonal appears complicated.

One subtlety in searching for local minima is that any given
potential function will usually not detect all possible families
of local minima that could occur for other potential functions.
For example, for five points in $\R^3$, the family with a square
opposite a point does not contain a local minimum for harmonic
energy.  One can attain a local minimum compared to the other
members of the family, but it will be a saddle point in the space
of all configurations. Nevertheless, the family does contain local
minima for some other completely monotonic potential functions
(such as $f(r)=1/r^s$ with $s$ large).

\subsection{Example}
\label{subsec:27in6}

\begin{table}
\begin{tabular}{ccccc}
Harmonic energy & Frequency & Parameters & Maximal cosine & Symmetries\\
\hline $111.0000000000$ & $99971504$ & $0$ & $0.2500000000$ & $51840$\\
$112.6145815185$ & $653$ & $9$ & $0.4306480635$ & $120$\\
$112.6420995468$ & $22993$ & $18$ & $0.3789599707$ & $24$\\
$112.7360209988$ & $10$ & $2$ & $0.4015602076$ & $1920$\\
$112.8896851626$ & $4840$ & $13$ & $0.4041651631$ & $48$\\
\\
\end{tabular}
\caption{Local minima for $27$ points in $\R^6$ (with frequencies
out of $10^8$ random trials).} \label{table:27}
\end{table}

For a concrete example, consider Table~\ref{table:27}, which shows
the results of $10^8$ random trials for $27$ points in $\R^6$ (all
decimal numbers in tables have been rounded). These parameters
were chosen because, as shown in \cite{CK1}, there is a unique
$27$-point universal optimum in $\R^6$, with harmonic energy
$111$; it is called the Schl\"afli configuration.  The column
labeled ``frequency'' tells how many times each local minimum
occurred. As one can see, the universal optimum occurred more than
99.97\% of the time, but we found a total of four others.

Strictly speaking, we have not proved that the local minima listed
in Table~\ref{table:27} (other than the Schl\"afli configuration)
even exist. They surely do, because we have computed them to five
hundred decimal places and checked that they are local minima by
numerically diagonalizing the Hessian matrix of the energy
function on the space of configurations. However, we used
high-precision floating point arithmetic, so this calculation does
not constitute a rigorous proof, although it leaves no reasonable
doubt. It is not at all clear whether there are additional local
minima. We have not found any, but the fact that one of the local
minima occurs only once in every ten million trials suggests that
there might be others with even smaller basins of attraction.

The local minimum with energy $112.736\dots$ stands out in two
respects besides its extreme rarity: it has many symmetries and it
depends on few parameters. That suggests that it should have a
simple description, and in fact it does, as a modification of the
universal optimum.  Only two inner products occur between distinct
points in the Schl\"afli configuration, namely $-1/2$ and $1/4$.
In particular, it is not antipodal, so one can define a new code
by replacing a single point $x$ with its antipode $-x$. The
remaining $26$ points can be divided into two clusters according
to their distances from $-x$. Immediately after replacing $x$ with
$-x$ the code will no longer be a local minimum, but if one allows
it to equilibrate a minimum is achieved.  (That is not obvious:
the code could equilibrate to a saddle point, because it is
starting from an unusual position.) All that changes is the
distances of the two clusters from $-x$, while the relative
positions within the clusters remain unchanged (aside from
rescaling). These two distances are the two parameters of the
code. The symmetries of the new code are exactly those of the
universal optimum that fix $x$, so the size of the symmetry group
is reduced by a factor of $27$.

The Schl\"afli configuration in $\R^6$ corresponds to the $27$
lines on a smooth cubic surface: there is a natural correspondence
between points in the configuration and lines on a cubic surface
so that the inner products of $-1/2$ occur between points
corresponding to intersecting lines. (This dates back to Schoutte
\cite{Sch}. See also the introduction to \cite{CK1} for a brief
summary of the correspondence.)  One way to view the other local
minima in Table~\ref{table:27} is as competitors to this classical
configuration.  It would be intriguing if they also had
interpretations or consequences in algebraic geometry, but we do
not know of any.

\section{Experimental phenomena}

\subsection{Analysis of Gram matrices}
\label{subsec:gramanalysis}

\begin{figure}
\begin{center}
\includegraphics{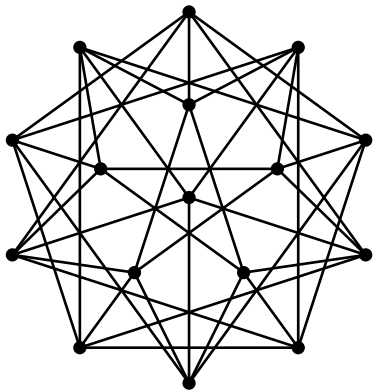}
\end{center}
\caption{The Clebsch graph.} \label{fig:clebsch}
\end{figure}

For an example of how one might analyze a Gram matrix, consider
the case of sixteen points in $\R^5$. This case also has a
universal optimum, in fact the smallest known one that is not a
regular polytope (although it is semiregular).  It is the
five-dimensional hemicube, which consists of half the vertices of
the cube.  More precisely, it contains the points $(\pm 1, \pm 1,
\pm 1, \pm 1, \pm 1)/\sqrt{5}$ with an even number of minus signs.
One can recover the full cube by including the antipode of each
point, so the symmetries of the five-dimensional hemicube consist
of half of those of the five-dimensional cube (namely, those that
preserve the hemicube, rather than exchanging it with its
complementary hemicube).

It is essentially an accident of five dimensions that the hemicube
is universally optimal.  Universal optimality also holds in lower
dimensions, but only because the hemicubes turn out to be familiar
codes (two antipodal points in two dimensions, a tetrahedron in
three dimensions, and a cross polytope in four dimensions).  In
six dimensions the hemicube appears to be an optimal spherical
code, but it does not minimize harmonic energy and is therefore
not universally optimal. In seven dimensions, and presumably all
higher dimensions, the hemicube is not even an optimal code.

The five-dimensional hemicube has the same structure as the
Clebsch graph (see Figure~\ref{fig:clebsch}).  The sixteen points
correspond to the vertices of the graph; two distinct points have
inner product $-3/5$ if they are connected by an edge in the graph
and $1/5$ otherwise.  This determines the Gram matrix and hence
the full configuration.

For the harmonic potential energy, the hemicube appears to be the
only local minimum with sixteen points in $\R^5$, but we do not
know how to prove that. To construct another local minimum, one
can attempt constructions such as moving a point to its antipode,
as in Subsection~\ref{subsec:27in6}, but they yield saddle points.
However, for other potential functions one sometimes finds other
local minima (we have found up to two other nontrivial local
minima). To illustrate the techniques from the previous section,
we will analyze one of them here.  It will turn out to have a
fairly simple conceptual description; our goal here is to explain
how to arrive at it, starting from the Gram matrix.

\begin{table}
\begin{center}
\begin{tabular}{c||ccc|ccc||ccc|ccc|ccc}
$1$ & $a$ & $a$ & $a$ & $a$ & $a$ & $a$ & $b$ & $b$ & $b$ & $b$ & $b$ & $b$ & $b$ & $b$ & $b$\\
\hline \hline $a$ & $1$ & $e$ & $e$ & $a^2$ & $a^2$ & $a^2$ & $d$ & $c$ & $c$ & $c$ & $c$ & $d$ & $c$ & $d$ & $c$\\
$a$ & $e$ & $1$ & $e$ & $a^2$ & $a^2$ & $a^2$ & $c$ & $d$ & $c$ & $d$ & $c$ & $c$ & $c$ & $c$ & $d$\\
$a$ & $e$ & $e$ & $1$ & $a^2$ & $a^2$ & $a^2$ & $c$ & $c$ & $d$ & $c$ & $d$ & $c$ & $d$ & $c$ & $c$\\
\hline $a$ & $a^2$ & $a^2$ & $a^2$ & $1$ & $e$ & $e$ & $d$ & $c$ & $c$ & $c$ & $d$ & $c$ & $c$ & $c$ & $d$\\
$a$ & $a^2$ & $a^2$ & $a^2$ & $e$ & $1$ & $e$ & $c$ & $d$ & $c$ & $c$ & $c$ & $d$ & $d$ & $c$ & $c$\\
$a$ & $a^2$ & $a^2$ & $a^2$ & $e$ & $e$ & $1$ & $c$ & $c$ & $d$ & $d$ & $c$ & $c$ & $c$ & $d$ & $c$\\
\hline\hline $b$ & $d$ & $c$ & $c$ & $d$ & $c$ & $c$ & $1$ & $f$ & $f$ & $f$ & $g$ & $g$ & $f$ & $g$ & $g$\\
$b$ & $c$ & $d$ & $c$ & $c$ & $d$ & $c$ & $f$ & $1$ & $f$ & $g$ & $f$ & $g$ & $g$ & $f$ & $g$\\
$b$ & $c$ & $c$ & $d$ & $c$ & $c$ & $d$ & $f$ & $f$ & $1$ & $g$ & $g$ & $f$ & $g$ & $g$ & $f$\\
\hline$b$ & $c$ & $d$ & $c$ & $c$ & $c$ & $d$ & $f$ & $g$ & $g$ & $1$ & $f$ & $f$ & $f$ & $g$ & $g$\\
$b$ & $c$ & $c$ & $d$ & $d$ & $c$ & $c$ & $g$ & $f$ & $g$ & $f$ & $1$ & $f$ & $g$ & $f$ & $g$\\
$b$ & $d$ & $c$ & $c$ & $c$ & $d$ & $c$ & $g$ & $g$ & $f$ & $f$ & $f$ & $1$ & $g$ & $g$ & $f$\\
\hline$b$ & $c$ & $c$ & $d$ & $c$ & $d$ & $c$ & $f$ & $g$ & $g$ & $f$ & $g$ & $g$ & $1$ & $f$ & $f$\\
$b$ & $d$ & $c$ & $c$ & $c$ & $c$ & $d$ & $g$ & $f$ & $g$ & $g$ & $f$ & $g$ & $f$ & $1$ & $f$\\
$b$ & $c$ & $d$ & $c$ & $d$ & $c$ & $c$ & $g$ & $g$ & $f$ & $g$ &
$g$ & $f$ & $f$ & $f$ & $1$
\end{tabular}
\end{center}
\caption{Gram matrix for $16$ points in $\R^5$; here $c = ab +
(1/2)\sqrt{(1-a^2)(1-b^2)/2}$, $d = ab - \sqrt{(1-a^2)(1-b^2)/2}$,
% The factor of (1/2) from c really should be missing in d; this is not a typo.
$e = (3a^2-1)/2$, $f = (3b^2-1)/2$, and $g = (3b^2+1)/4$.}
\label{table:gram16in5}
\end{table}

The specific example we will analyze arises as a local minimum for
the potential function $r \mapsto (4-r)^{12}$.  It is specified by
Table~\ref{table:gram16in5} with
$$
a \approx -0.499890010934
$$
and
$$
b \approx 0.201039702365
$$
(the lines in the table are just for visual clarity).

The first step is to recognize the structure in the Gram matrix.
Table~\ref{table:gram16in5} highlights this structure, but of
course it takes effort to bring the Gram matrix into such a simple
form (by recognizing algebraic relations between the Gram matrix
entries and reordering the points so as to emphasize patterns).
The final form of the Gram matrix exhibits the configuration as
belonging to a family specified by two parameters $a$ and $b$ with
absolute value less than $1$.  As described in the table's
caption, all the other inner products are simple algebraic
functions of $a$ and $b$.  To check that this Gram matrix
corresponds to an actual code in $S^4$, it suffices to verify that
its eigenvalues are $0$ ($11$ times), $1+6a^2+9b^2$, and
$(15-(6a^2+9b^2))/4$ ($4$ times): there are only five nonzero
eigenvalues and they are clearly positive.

Table~\ref{table:gram16in5} provides a complete description of the
configuration, but it is unilluminating.  To describe the code using
elegant coordinates, one must have a more conceptual understanding of
it.  A first step in that direction is the observation that the first
point in Table~\ref{table:gram16in5} has inner product $a$ or $b$ with
every other point.  In other words, the remaining $15$ points lie on
two parallel four-dimensional hyperplanes, equidistant from the first
point. A natural guess is that as $a$ and $b$ vary, the structures
within these hyperplanes are simply rescaled as the corresponding cross
sections of the sphere change in size, and some calculation verifies
that this guess is correct.

To understand these two structures and how they relate to each
other, set $a=b=0$ so that they form a $15$-point configuration in
$\R^4$.  Its Gram matrix is of course obtained by removing the
first row and column of Table~\ref{table:gram16in5} and setting
$a=b=0$, $e=f=-1/2$, $g=1/4$, $c=\sqrt{2}/4$, and $d=-\sqrt{2}/2$.
The two substructures consist of the first six points and the last
nine, among the fifteen remaining points.

Understanding the $16$-point codes in $\R^5$ therefore simply
comes down to understanding this single $15$-point code in $\R^4$.
(It is also the $15$-point code from Table~\ref{table:balanced}.
Incidentally, Sloane's tables \cite{S} show that it is not an
optimal spherical code.) The key to understanding it is choosing
the right coordinates. The first six points form two orthogonal
triangles, and they are the simplest part of this configuration,
so it is natural to start with them.

Suppose the points $v_1,v_2,v_3$ and $v_4,v_5,v_6$ form two
orthogonal equilateral triangles in a four-dimensional vector
space. The most natural coordinates to choose for the vector space
are the inner products with these six points.  Of course the sum
of the three inner products with any triangle must vanish (because
$v_1+v_2+v_3=v_4+v_5+v_6=0$), so there are only four independent
coordinates, but we prefer not to break the symmetry by discarding
two coordinates.

The other nine points in the configuration are determined by their
inner products with $v_1,\dots,v_6$.  Each of them will have inner
product $d$ with one point in each triangle and $c$ with the
remaining two points. As pointed out above we must have $d+2c=0$,
and in fact $d=-\sqrt{2}/2$ and $c=\sqrt{2}/4$ because the points
are all unit vectors.  Note that one can read off all this
information from the $c$ and $d$ entries in
Table~\ref{table:gram16in5}.

There is an important conceptual point in the last part of this
analysis.  Instead of focusing on the internal structure among the
last nine points, it is most fruitful to study how they relate to
the previously understood subconfiguration of six points. However,
once one has a complete description, it is important to examine
the internal structure as well.

The pattern of connections among the last nine points in
Table~\ref{table:gram16in5} is described by the Paley graph on
nine vertices, which is the unique strongly regular graph with
parameters $(9,4,1,2)$.  (The Paley graph is isomorphic to its own
complement, so the edges could correspond to inner product either
$f$ or $g$.) Strongly regular graphs, and more generally
association schemes, frequently occur as substructures of
minimal-energy configurations.  It is remarkable to see such
highly ordered structures spontaneously occurring via energy
minimization.

\subsection{Other small examples}

To illustrate some of the other phenomena that can occur, in this
subsection we will analyze the case of $12$ points in $\R^4$.  We
have observed two families of local minima, both of which are
slightly more subtle than the previous examples.

For $0 < a < 1/2$, set $b = a - 1$, $c = -3a + 1$, and $d = 4a -
1$, and consider the Gram matrix shown in
Table~\ref{table:gram12in4a} (its nonzero eigenvalues are $12a$
and $6-12a$, each with multiplicity $2$). Unlike the example in
Subsection~\ref{subsec:gramanalysis}, the symmetry group acts
transitively on the points, so there are no distinguished points
to play a special role in the analysis. Nevertheless, one can
analyze it as follows.

\begin{table}
\begin{center}
\begin{tabular}{ccc|ccc|ccc|ccc}
$1$ & $c$ & $c$ & $d$ & $b$ & $b$ & $-2a$ & $a$ & $a$ & $-2a$ & $a$ & $a$\\
$c$ & $1$ & $c$ & $b$ & $d$ & $b$ & $a$ & $-2a$ & $a$ & $a$ & $-2a$ & $a$\\
$c$ & $c$ & $1$ & $b$ & $b$ & $d$ & $a$ & $a$ & $-2a$ & $a$ & $a$ & $-2a$\\
\hline$d$ & $b$ & $b$ & $1$ & $c$ & $c$ & $-2a$ & $a$ & $a$ & $-2a$ & $a$ & $a$\\
$b$ & $d$ & $b$ & $c$ & $1$ & $c$ & $a$ & $-2a$ & $a$ & $a$ & $-2a$ & $a$\\
$b$ & $b$ & $d$ & $c$ & $c$ & $1$ & $a$ & $a$ & $-2a$ & $a$ & $a$ & $-2a$\\
\hline$-2a$ & $a$ & $a$ & $-2a$ & $a$ & $a$ & $1$ & $c$ & $c$ & $d$ & $b$ & $b$\\
$a$ & $-2a$ & $a$ & $a$ & $-2a$ & $a$ & $c$ & $1$ & $c$ & $b$ & $d$ & $b$\\
$a$ & $a$ & $-2a$ & $a$ & $a$ & $-2a$ & $c$ & $c$ & $1$ & $b$ & $b$ & $d$\\
\hline$-2a$ & $a$ & $a$ & $-2a$ & $a$ & $a$ & $d$ & $b$ & $b$ & $1$ & $c$ & $c$\\
$a$ & $-2a$ & $a$ & $a$ & $-2a$ & $a$ & $b$ & $d$ & $b$ & $c$ & $1$ & $c$\\
$a$ & $a$ & $-2a$ & $a$ & $a$ & $-2a$ & $b$ & $b$ & $d$ & $c$ &
$c$ & $1$
\end{tabular}
\end{center}
\caption{Gram matrix for $12$ points in $\R^4$; here $0 < a <
1/2$, $b = a - 1$, $c = -3a + 1$, and $d = 4a - 1$.}
\label{table:gram12in4a}
\end{table}

Let $v_1,v_2,v_3 \in S^1$ be the vertices of an equilateral
triangle in $\R^2$, and let $v_4$ and $v_5$ be unit vectors that
are orthogonal to each other and to each of $v_1$, $v_2$, and
$v_3$. For $0 < \alpha < 1$, consider the twelve points $\alpha
v_i \pm \sqrt{1-\alpha^2} v_4$ and $-\alpha v_i \pm
\sqrt{1-\alpha^2} v_5$ with $1 \le i \le 3$.  If one sets $a =
\alpha^2/2$ then they have Table~\ref{table:gram12in4a} as a Gram
matrix.

The Gram matrix shown in Table~\ref{table:gram12in4b} is quite
different.  There, $0 < a < 1/3$, $b = 1-12a^2$, $c = 6a^2-1$, and
$d = 18a^2-1$.  The nonzero eigenvalues are $4/3+24a^2$ (with
multiplicity $3$) and $8-72a^2$, which are positive because $a <
1/3$.

\begin{table}
\begin{center}
\begin{tabular}{cccc|cccc|cccc}
$1$ & $-1/3$ & $-1/3$ & $-1/3$ & $-3a$ & $a$ & $a$ & $a$ & $-3a$ & $a$ & $a$ & $a$\\
$-1/3$ & $1$ & $-1/3$ & $-1/3$ & $a$ & $-3a$ & $a$ & $a$ & $a$ & $-3a$ & $a$ & $a$\\
$-1/3$ & $-1/3$ & $1$ & $-1/3$ & $a$ & $a$ & $-3a$ & $a$ & $a$ & $a$ & $-3a$ & $a$\\
$-1/3$ & $-1/3$ & $-1/3$ & $1$ & $a$ & $a$ & $a$ & $-3a$ & $a$ & $a$ & $a$ & $-3a$\\
\hline
$-3a$ & $a$ & $a$ & $a$ & $1$ & $b$ & $b$ & $b$ & $d$ & $c$ & $c$ & $c$\\
$a$ & $-3a$ & $a$ & $a$ & $b$ & $1$ & $b$ & $b$ & $c$ & $d$ & $c$ & $c$\\
$a$ & $a$ & $-3a$ & $a$ & $b$ & $b$ & $1$ & $b$ & $c$ & $c$ & $d$ & $c$\\
$a$ & $a$ & $a$ & $-3a$ & $b$ & $b$ & $b$ & $1$ & $c$ & $c$ & $c$ & $d$\\
\hline
$-3a$ & $a$ & $a$ & $a$ & $d$ & $c$ & $c$ & $c$ & $1$ & $b$ & $b$ & $b$\\
$a$ & $-3a$ & $a$ & $a$ & $c$ & $d$ & $c$ & $c$ & $b$ & $1$ & $b$ & $b$\\
$a$ & $a$ & $-3a$ & $a$ & $c$ & $c$ & $d$ & $c$ & $b$ & $b$ & $1$ & $b$\\
$a$ & $a$ & $a$ & $-3a$ & $c$ & $c$ & $c$ & $d$ & $b$ & $b$ & $b$
& $1$
\end{tabular}
\end{center}
\caption{Gram matrix for $12$ points in $\R^4$; here $0 < a <
1/3$, $b = 1-12a^2$, $c = 6a^2-1$, and $d = 18a^2-1$.}
\label{table:gram12in4b}
\end{table}

In this Gram matrix the first four points form a distinguished
tetrahedron, and the remaining eight points form two identical
tetrahedra. They lie in hyperplanes parallel to and equidistant
from the (equatorial) hyperplane containing the distinguished
tetrahedron. If one sets $a=1/3$, then all three tetrahedra lie in
the same hyperplane, with $b=-1/3$, $c=-1/3$, $d=1$, and $-3a=-1$.
In particular, one can see that the two parallel tetrahedra are in
dual position to the distinguished tetrahedron. As the parameter
$a$ varies, all that changes is the distance between the parallel
hyperplanes. (As $a$ tends to zero some points coincide.  One
could also use $a$ between $0$ and $-1/3$, but that corresponds to
using parallel tetrahedra oriented the same way, instead of
dually, which generally yields higher potential energy.)

This sort of layered structure occurs surprisingly often.  One
striking example is $74$ points in $\R^5$.  The best such
spherical code known consists of a regular $24$-cell on the
equatorial hyperplane together with two dual $24$-cells on
parallel hyperplanes as well as the north and south poles.  If one
chooses the two parallel hyperplanes to have inner products $\pm
\sqrt{\sqrt{5}-2}$ with the poles, then the cosine of the minimal
angle is exactly $(\sqrt{5}-1)/2$.  That agrees numerically with
Sloane's tables \cite{S} of the best codes known, but of course
there is no proof that it is optimal.

There is almost certainly no universally optimal $12$-point
configuration in $\R^4$.  Aside from some trivial examples for
degenerate potential functions, the two cases we have analyzed in
this subsection are the only two types of local minima we have
observed.  For $f(r) = (4-r)^k$ with $k \in \{1,2\}$ they both
achieve the same minimal energy (along with a positive-dimensional
family of other configurations).  For $3 \le k \le 9$ the first
family appears to achieve the global minimum, while for $k \ge 10$
the second appears to.  As $k$ tends to infinity the energy
minimization problem turns into the problem of finding the optimal
spherical code.  That problem appears to be solved by taking
$a=1/4$ in the second family, so that the minimal angle has cosine
$1/4$, which agrees with Sloane's tables \cite{S}.

We conjecture that one or the other of these two families
minimizes each completely monotonic potential function.  This
conjecture is somewhat difficult to test, but we are not aware of
any counterexamples.

The examples we have analyzed so far illustrate three basic
principles:
\begin{enumerate}
\item Small or medium-sized local minima tend to occur
in low-dimensional families as one varies the potential function.
The dimension is not usually as low as in these examples, but it
is typically far lower than the dimension of the space of all
configurations (see Figure~\ref{fig:params}).

\item These families frequently contain surprisingly symmetrical
substructures (such as regular polytopes or configurations
described by strongly regular graphs or other association
schemes).

\item The same substructures and construction methods occur in many different
families.
\end{enumerate}

\subsection{$2n+1$ points in $\R^n$}

Optimal spherical codes are known for up to $2n$ points in $\R^n$
(see Theorem~6.2.1 in \cite{B}), but not for $2n+1$ points, except
in $\R^2$ and $\R^3$.  Here we present a natural conjecture for
all dimensions.

These codes consist of a single point we call the north pole
together with two $n$-point simplices on hyperplanes parallel to
the equator; the simplices are in dual position relative to each
other.  Each point in the simplex closer to the north pole will
have inner product $\alpha$ with the north pole, and the inner
product between any two points in the further simplex will be
$\alpha$.  The number $\alpha$ can be chosen so that each point in
either one of the simplices has inner product $\alpha$ with each
point in the other simplex except the point furthest from it. To
achieve that, $\alpha$ must be the unique root between $0$ and
$1/n$ of the cubic equation
$$
(n^3-4n^2+4n)x^3 - n^2x^2 -nx + 1 = 0.
$$
As $n \to \infty$, $\alpha = 1/n - \sqrt{2}/n^{3/2} + O(1/n^2)$.

Let $\mathcal{C}_n \subset S^{n-1}$ be this spherical code, with
$\alpha$ chosen as above. The cosine of the minimal angle in
$\mathcal{C}_n$ is $\alpha$.

\begin{conjecture}
\label{conj:opt} For each $n \ge 2$, the code $\mathcal{C}_n$ is
an optimal spherical code.  Furthermore, every optimal
$(2n+1)$-point code in $S^{n-1}$ is isometric to $\mathcal{C}_n$.
\end{conjecture}

On philosophical grounds it seems reasonable to expect to be able
to prove this conjecture: most of the difficulty in packing
problems comes from the idiosyncrasies of particular spaces and
dimensions, so when a phenomenon occurs systematically one expects
a conceptual reason for it.  However, we have made no serious
progress towards a proof.

One can also construct $\mathcal{C}_n$ as follows. Imagine adding
one point to a regular cross polytope by placing it in the center
of a facet. The vertices of that facet form a simplex equidistant
from the new point, as do the vertices of the opposite facet.  The
structure is identical to the code $\mathcal{C}_n$, except for the
distances from the new point, and the proper distances can be
obtained by allowing the code to equilibrate with respect to
increasingly steep potential functions.

It appears that for $n>2$ these codes do not minimize harmonic
energy, so they are not universally optimal. When $n=4$, something
remarkable occurs with the (conjectured) minimum for harmonic
energy.  That configuration consists of a regular pentagon
together with two pairs of antipodal points that are orthogonal to
each other and the pentagon.  If one uses gradient descent to
minimize harmonic energy, it seems to converge with probability
$1$ to this configuration, but the convergence is very slow, much
slower than for any other harmonic energy minimum we have found.
The reason is that this configuration is a degenerate minimum for
the harmonic energy, in the sense that the Hessian matrix has more
zero eigenvalues than one would expect.

Each of the nine points has three degrees of freedom, so the
Hessian matrix has twenty-seven eigenvalues. Specifically, they
are $0$ (ten times), $4$, $7/4$ (twice), $9/2$ (four times), $9$
(twice), $25/8\pm\sqrt{209}/8$ (twice), and $31/8 \pm
\sqrt{161}/8$ (twice).  Six of the zero eigenvalues are
unsurprising, because they come from the problem's invariance
under the six-dimensional Lie group $O(4)$, but the remaining four
are surprising indeed.

The corresponding eigenvectors are infinitesimal displacements of
the nine points that produce only a fourth-order change in energy,
rather than the expected second-order change. To construct them,
do not move the antipodal pairs of points at all, and move the
pentagon points orthogonally to the plane of the pentagon.  Each
must be displaced by $(1-\sqrt{5})/2$ times the sum of the
displacements of its two neighbors.  This yields a
four-dimensional space of displacements, which are the surprising
eigenvectors.

This example is noteworthy because it shows that harmonic energy
is not always a Morse function on the space of all configurations.
One might hope to apply Morse theory to understand the
relationship between critical points for energy and the topology
of the configuration space, but the existence of degenerate
critical points could substantially complicate this approach.

\subsection{$2n+2$ points in $\R^n$}
\label{subsec:diplo}

After seeing a conjecture for the optimal $(2n+1)$-point code in
$S^{n-1}$, it is natural to wonder about $2n+2$ points.  A first
guess is the union of a simplex and its dual simplex (in other
words, the antipodal simplex), which was named the diplo-simplex
by Conway and Sloane \cite{CS2}. One can prove using the linear
programming bounds for real projective space that this code is the
unique optimal antipodal spherical code of its size and dimension
(see Chapter~9 of \cite{CS}), but for $n>2$ it is not even locally
optimal as a spherical code (see Appendix~\ref{appendix:diplo})
and we do not have a conjecture for the true answer.

For the problem of minimizing harmonic energy, the diplo-simplex
is suboptimal for $3 \le n \le 5$ but appears optimal for all
other $n$.

\begin{figure}
\begin{center}
\includegraphics{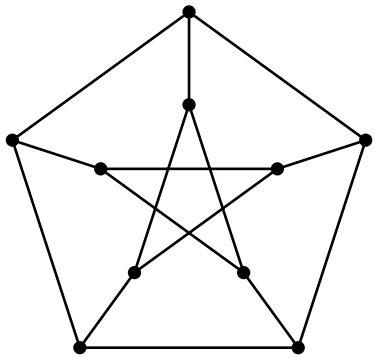}
\end{center}
\caption{The Petersen graph.} \label{fig:petersen}
\end{figure}

One particularly elegant case is when $n=4$.  The midpoints of the
edges of a regular simplex form a $10$-point code in $S^3$ with maximal
inner product $1/6$, and Bachoc and Vallentin \cite{BV} have proved
that it is the unique optimal spherical code. It is also the kissing
configuration of the five-dimensional hemicube (the universally optimal
$16$-point configuration in $\R^5$). In other words, it consists of the
ten nearest neighbors of any point in that code. This code appears to
minimize harmonic energy, but it is not the unique minimum: two
orthogonal regular pentagons have the same harmonic energy.

As pointed out in the introduction of \cite{CK1}, this code is not
universally optimal, but it nevertheless seems to be an
exceedingly interesting configuration.  Only the inner products
$-2/3$ and $1/6$ occur (besides $1$, of course). If one forms a
graph whose vertices are the points in the code and whose edges
correspond to pairs of points with inner product $-2/3$, then the
result is the famous Petersen graph (Figure~\ref{fig:petersen}).

Like the nontrivial universal optima in dimensions $5$ through $8$,
this code consists of the vertices of a semiregular polytope that has
simplices and cross polytopes as facets, with a simplex and two cross
polytopes meeting at each face of codimension~$3$. Its kissing
configuration is also semiregular, with square and triangular facets,
but it is a suboptimal code (specifically, a triangular prism).

\subsection{$48$ points in $\R^4$}
\label{subsec:48in4}

One of the most beautiful configurations we have found is a
$48$-point code in $\R^4$. The points form six octagons that map
to the vertices of a regular octahedron under the Hopf map from
$S^3$ to $S^2$.  Recall that if we identify $\R^4$ with $\C^2$
using the inner product $\langle x,y \rangle =
\mathop{\textup{Re}} \bar x ^t y$ on $\C^2$, then the Hopf map
sends $(z,w)$ to $z/w \in \C \cup \{\infty\}$, which we can
identify with $S^2$ via stereographic projection to a unit sphere
centered at the origin.  The fibers of the Hopf map are the
circles given by intersecting $S^3$ with the complex lines in
$\C^2$.

Sloane, Hardin, and Cara \cite{SHC} found a spherical $7$-design of
this form, consisting of two dual $24$-cells, and it has the same
minimal angle as our code (which is the minimal angle in an octagon),
but it is a different code. In $\C^2$, the Sloane-Hardin-Cara code is
the union of the orbits under multiplication by eighth roots of unity
of the points $(1,0)$, $(0,1)$, $(\pm 1, 1)/\sqrt{2}$, and $(\pm i,
1)/\sqrt{2}$.  Our code is the union of the orbits of $(1,0)$, $(0,1)$,
$(\pm\zeta,\zeta)/\sqrt{2}$, and $(\pm i\zeta^2,\zeta^2)/\sqrt{2}$,
where $\zeta = e^{\pi i/12}$. Each octagon has been rotated by a
multiple of $\pi/12$ radians.  Because a regular octagon is invariant
under rotation by $\pi/4$ radians, there are only three distinct
rotations by multiples of $\pi/12$.  Each such rotation occurs for the
octagons lying over two antipodal vertices of the octahedron in the
base space $S^2$ of the Hopf fibration.

It is already remarkable that performing these rotations yields a
balanced configuration with lower harmonic energy than the union
of the $24$-cell and its dual, but the structure of the code's
convex hull is especially noteworthy.  The facets can be computed
using the program Polymake \cite{GJ}. The facets of the dual
$24$-cell configuration are $288$ irregular tetrahedra, all
equivalent under the action of the symmetry group (and each
possessing $8$ symmetries).  By contrast, our code has $128$
facets forming two orbits under the symmetry group: one orbit of
$96$ irregular tetrahedra and one of $32$ irregular octahedra. The
irregular octahedra are obtained from regular ones by rotating one
of the facets, which are equilateral triangles, by an angle of
$\pi/12$. We will use the term ``twisted facets'' to denote the
rotated facet and its opposite facet (by symmetry, either one
could be viewed as rotated relative to the other).

The octahedra in our configuration meet other octahedra along
their twisted facets and simplices along their other facets.
Grouping the octahedra according to adjacency therefore yields
twisted chains of octahedra. Each chain consists of eight
octahedra, and they span the $3$-sphere along great circles.  The
total twist amounts to $8\pi/12=2\pi/3$, from which it follows
that the chains close with facets aligned correctly.  The $32$
octahedra form four such chains, and the corresponding great
circles are fibers in the same Hopf fibration as the vertices of
the configuration.  These Hopf fibers map to the vertices of a
regular tetrahedron in $S^2$.  It is inscribed in the cube dual to
the octahedron formed by the images of the vertices of the code.

Another way to view the facets of this polytope, or any spherical
polytope, is as holes in the spherical code.  More precisely, the
(outer) facet normals of any full-dimensional polytope inscribed
in a sphere are the holes in the spherical code (i.e., the points
on the sphere that are local maxima for distance from the code).
The normals of the octahedral facets are the deep holes in this
code (i.e., the points at which the distance is globally
maximized). Notice that these points are defined using the
intrinsic geometry of the sphere, rather than relying on its
embedding in Euclidean space.

The octahedral facets of our code can be thought of as more
important than the tetrahedral facets.  The octahedra appear to us
to have prettier, clearer structure, and once they have been
placed, the entire code is determined (the tetrahedra simply fill
the gaps). This idea is not mathematically precise, but it is a
common theme in many of our calculations: when we examine the
facet structure of a balanced code, we often find a small number
of important facets and a large number of less meaningful ones.

\subsection{Hopf structure}
\label{subsec:Hopf}

As in the previous example, many notable codes in $S^3$, $S^7$, or
$S^{15}$ can be understood using the complex, quaternionic, or
octonionic Hopf maps (see for example \cite{D} and \cite{AP-G3}).
In this subsection, we describe this phenomenon for the regular
$120$-cell and $600$-cell in $S^3$.  The Hopf structure on the
$600$-cell is mathematical folklore, but we have not been able to
locate it in the published literature, while the case of the
$120$-cell is more subtle and may not have been previously
examined.

\begin{figure}
\begin{center}
\includegraphics{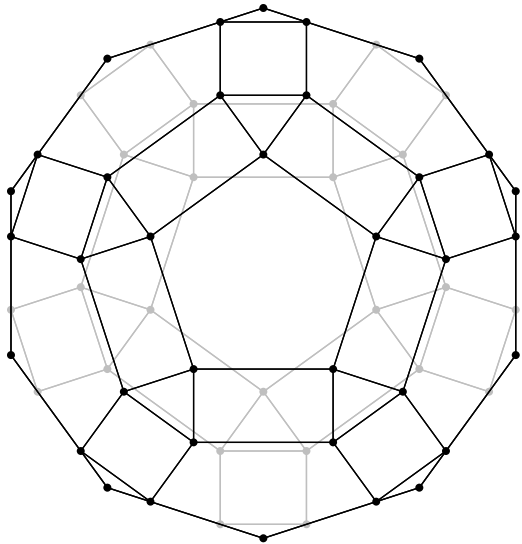}
\end{center}
\caption{The $60$-point polytope in $\R^3$ over which the regular
$120$-cell fibers. } \label{fig:poly60}
\end{figure}

The $H_4$ reflection group (which is the symmetry group of both
polytopes) contains elements of order $10$ that act on $\R^4$ with no
fixed points other than the origin.  If one chooses such an element,
then $\R^4$ has the structure of a two-dimensional complex vector space
such that this element acts via multiplication by a primitive $10$-th
root of unity.  The orbits are regular $10$-gons lying in Hopf fibers.
In the case of the regular $600$-cell, this partitions the $120$
vertices into $12$ regular $10$-gons lying in Hopf fibers over the
vertices of a regular icosahedron in $S^2$. For the regular $120$-cell
(with $600$ vertices), the corresponding polyhedron in $S^2$ has $60$
vertices, but it is far from obvious what it is. We know of no way to
determine it without calculation, but computing with coordinates
reveals that it is a distorted rhombicosidodecahedron, with the square
facets replaced by golden rectangles. Specifically, its facets are $12$
regular pentagons, $20$ equilateral triangles, and $30$ golden
rectangles. The golden rectangles meet pentagons along their long edges
and triangles along their short edges. Figure~\ref{fig:poly60} shows
the orthogonal projection into the plane containing a pentagonal face
(gray vertices and edges are on the far side of the polyhedron).

\subsection{Facet structure of universal optima}
\label{subsec:facet}

The known low-dimensional universal optima (through dimension $8$)
are all regular or semiregular polytopes, whose facets are well
known to be regular simplices or cross polytopes. However, there
seems to have been little investigation of the facets of the
higher-dimensional universal optima from
Table~\ref{table:universal}.  In this subsection we will look at
the smallest higher-dimensional cases: $\universal{100}{22}$,
$\universal{112}{21}$, $\universal{162}{21}$,
$\universal{275}{22}$, $\universal{552}{23}$, and
$\universal{891}{22}$ (recall that $\universal{N}{n}$ denotes the
$N$-point code in $\R^n$ from Table~\ref{table:universal}, when it
is unique).  Each of the first four is a two-distance set given by
a spectral embedding of a strongly regular graph. (Recall that a
spectral embedding is obtained by orthogonally projecting the
standard orthonormal basis into an eigenspace of the adjacency
matrix of the graph.)  The last two have three distances between
distinct points.

These codes have enormous numbers of facets (more than
seventy-five trillion for $\universal{552}{23}$), so it is not
feasible to find the facets using general-purpose methods.
Instead, one must make full use of the large symmetry groups of
these configurations. With Dutour Sikiri\' c's package Polyhedral
\cite{DS} for the program GAP \cite{GAP}, that can be done for
these configurations. We have used it to compute complete lists of
orbits of facets under the action of the symmetry group. (The
results are rigorous, because we use exact coordinates for the
codes.  In particular, when necessary we use the columns of the
Gram matrices to embed scalar multiples of these codes
isometrically into high-dimensional spaces using only rational
coordinates.) Of course, the results of this computation then
require analysis by hand to reveal their structure.

For an introductory example, it is useful to review the case of the
five-dimensional hemicube (see Subsection~\ref{subsec:gramanalysis}).
It has ten obvious facets contained in the ten facets of the cube. Each
is a four-dimensional hemicube, i.e., a regular cross polytope.  The
remaining facets are regular simplices (one opposite each point in the
hemicube).

One can view the five-dimensional hemicube as an antiprism formed by
two four-dimensional cross polytopes in parallel hyperplanes. The cross
polytopes are arranged so that each one's vertices point towards deep
holes of the other.  (The deep holes of a cross polytope are the
vertices of the dual cube, and in four dimensions the vertices of that
cube consist of two cross polytopes. The fact that the deep holes of a
four-dimensional cross polytope contain another such cross polytope is
crucial for this construction to make sense.) Of course, the distance
between the parallel hyperplanes is chosen so as to maximize the
minimal distance. What is remarkable about this antiprism is that it is
far more symmetrical than one might expect: normally the two starting
facets of an antiprism play a very different role from the facets
formed when taking the convex hull, but in this case extra symmetries
occur. The simplest case of such extra symmetries is the construction
of a cross polytope as an antiprism made from two regular simplices in
dual position.

The three universal optima $\universal{100}{22}$,
$\universal{112}{21}$, and $\universal{162}{21}$ are each given by
an unusually symmetric antiprism construction analogous to that of
the hemicube. In each case, the largest facets (i.e., those
containing the most vertices) contain half the vertices. These
facets are themselves spectral embeddings of strongly regular
graphs (the Hoffman-Singleton graph, the Gewirtz graph, and the
unique $(81,20,1,6)$ strongly regular graph). Within the universal
optima, the largest facets occur in pairs in parallel hyperplanes,
and the vertices of each facet in a pair point towards holes in
the other.  These holes belong to a single orbit under the
symmetry group of the facet, and that orbit is the disjoint union
of several copies of the vertices of the facet: two copies for the
Hoffman-Singleton and Gewirtz cases and four in the third case.
These holes are the deepest holes in the Hoffman-Singleton case;
in the other two cases, they are not quite the deepest holes
(there are not enough deep holes for the construction to work
using them).

Brouwer and Haemers \cite{BH1,BH2} discovered the underlying
combinatorics of these constructions (i.e., that the strongly
regular graphs corresponding to the universal optima can be
naturally partitioned into two identical graphs).  However, the
geometric interpretation as antiprisms appears to be new.

The universal optima $\universal{100}{22}$, $\universal{112}{21}$, and
$\universal{162}{21}$ are antiprisms, but that cannot possibly be true
for $\universal{275}{22}$, because $275$ is odd. Instead, the
McLaughlin configuration $\universal{275}{22}$ is analogous to the
Schl\"afli configuration $\universal{27}{6}$. Both are two-distance
sets.  In the Schl\"afli configuration, the neighbors of each point
form a five-dimensional hemicube and the non-neighbors form a
five-dimensional cross polytope.  Both the hemicube and the cross
polytope are unusually symmetric antiprisms, and their vertices point
towards each other's deep holes.  (The deep holes of the hemicube form
a cross polytope, and those of the cross polytope form a cube
consisting of two hemicubes.)  The McLaughlin configuration is
completely analogous: the neighbors of each point form
$\universal{162}{21}$ and the non-neighbors form $\universal{112}{21}$.
They point towards each other's deep holes; this is possible because
the deep holes of $\universal{112}{21}$ consist of four copies of
$\universal{162}{21}$, and its deep holes consist of two copies of
$\universal{112}{21}$. Furthermore, the deep holes in these two
universal optima are of exactly the same depth (i.e., distance to the
nearest point in the code), as is also the case for the
five-dimensional cross polytope and hemicube used to form the
Schl\"afli configuration.

The Schl\"afli and McLaughlin configurations both have the
property that their deep holes are the antipodes of their
vertices.  Thus, it is natural to form antiprisms from two
parallel copies of them, with vertices pointed at each other's
deep holes. That yields antipodal configurations of $54$ points in
$\R^7$ and $550$ points in $\R^{23}$.  If one also includes the
two points orthogonal to the parallel hyperplanes containing the
original two copies, then this construction gives the universal
optima $\universal{56}{7}$ and $\universal{552}{23}$.

\begin{table}
\begin{center}
\begin{tabular}{cc|cc}
Vertices & Number of orbits & Vertices & Number of orbits\\
\hline
$22$ & $92$ & $30$ & $1$\\
$23$ & $13$ & $31$ & $1$\\
$24$ & $6$ & $36$ & $1$\\
$25$ & $3$ & $42$ & $1$\\
$27$ & $3$ & $50$ & $1$\\
$28$ & $1$\\
\\ % Bad to skip line?
\end{tabular}
\end{center}
\caption{Number of orbits of facets of different sizes in the
Higman-Sims configuration
$\universal{100}{22}$.}\label{table:22-100facets}
\end{table}

Each high-dimensional universal optimum has many types of facets
of different sizes.  For example, the facets of the Higman-Sims
configuration $\universal{100}{22}$ form $123$ orbits under the
action of the symmetry group (see Table~\ref{table:22-100facets}).
The largest facets, which come from the Hoffman-Singleton graph as
described above, are by far the most important, but each type of
facet appears to be of interest.  They are often more subtle than
one might expect. For example, it is natural to guess that the
facets with $42$ vertices would be regular cross polytopes, based
on the number of vertices, but they are not.  Instead, when
rescaled to the unit sphere they have the following structure:

The facets with $42$ vertices are two-distance sets on the unit
sphere in $\R^{21}$, with inner products $1/29$ and $-13/29$.  If
we define a graph on the vertices by letting edges correspond to
pairs with inner product $-13/29$, then this graph is the
bipartite incidence graph for points and lines in the projective
plane $\Proj^2(\F_4)$.  To embed this graph in $\R^{21}$,
represent the $21$ points in $\Proj^2(\F_4)$ as the permutations
of $(a,b,\dots,b)$, where $a^2+20b^2=1$ and $2ab+19b^2=1/29$.
Specifically, take $a = 0.9977\dots$ and $b = 0.0151\dots$ (these
are fourth degree algebraic numbers). Choose $c$ and $d$ so that
$5c^2+16d^2=1$ and $8cd+c^2+12d^2=1/29$ (specifically, take $c =
-0.4362\dots$ and $d = 0.0550\dots$).  Then embed the $21$ lines
into $\R^{21}$ as permutations of $(c,c,c,c,c,d,\dots,d)$, where
the five $c$ entries correspond to the points contained in the
line. This embedding gives the inner products of $1/29$ and
$-13/29$, as desired (and in fact those are the only inner
products for which a construction of this form is possible).

As shown in Table~\ref{table:22-100facets}, there are $92$
different types of simplicial facets in the Higman-Sims
configuration.  One orbit consists of regular simplices: for each
point in the configuration, the $22$ points at the furthest
distance from it form a regular simplex.  All the other simplices
are irregular. Nine orbits consist of simplices with no symmetries
whatsoever, and the remaining ones have some symmetries but not
the full symmetric group.

The universal optima $\universal{552}{23}$ and
$\universal{891}{22}$ have more elaborate facet structures, but we
have completely classified their facets (which form $116$ and
$422$ orbits, respectively). The facets corresponding to their
deep holes form single orbits, consisting of $\universal{100}{22}$
in the first case and $\universal{162}{21}$ in the second.  These
results can all be understood in terms of the standard embeddings
of these configurations into the Leech lattice, as follows:

Let $v$ be any vector with norm $6$ in the Leech lattice
$\Lambda_{24}$.  Among the $196560$ minimal vectors in
$\Lambda_{24}$ (those with norm $4$), there are $552$ minimal
vectors $w$ satisfying $|w-v|^2=4$, and they form a copy of the
$552$-point universal optimum.  This shows that
$\universal{552}{23}$ is a facet of $\universal{196560}{24}$.
Taking kissing configurations shows that $\universal{275}{22}$ is
a facet of $\universal{4600}{23}$ and that $\universal{162}{21}$
is a facet of $\universal{891}{22}$.  We conjecture that each of
these facets corresponds to a deep hole in the code, and that all
of the deep holes arise this way, but we have not proved this
conjecture beyond $\universal{891}{22}$.  The
$\universal{100}{22}$ facets of $\universal{552}{23}$ can also be
seen in this picture: given two vectors $v_1,v_2 \in \Lambda_{24}$
with $|v_1|^2=|v_2|^2=6$ and $|v_1-v_2|^2=4$, the corresponding
$\universal{552}{23}$ facets of $\universal{196560}{24}$ intersect
in a $\universal{100}{22}$ facet of $\universal{552}{23}$, which
corresponds to a deep hole.

\subsection{$96$ points in $\R^9$}
\label{subsec:96in9}

Another intriguing code that arose in our computer searches is a
$96$-point code in $\R^9$ (see Table~\ref{table:balanced}).  This
code was known previously: it is mentioned but not described in
Table~9.2 of \cite{CS}, which refers to a paper in preparation
that never appeared, and it is described in Appendix~D of
\cite{EZ}.  Here we describe it in detail, with a different
approach from that in \cite{EZ}.

The code is not universally optimal, but it is balanced and it
appears to be an optimal spherical code. What makes it noteworthy
is that the cosine of its minimal angle is $1/3$.  Any such code
corresponds to an arrangement of unit balls in $\R^{10}$ that are
all tangent to two fixed, tangent balls, where the interiors of
the balls are not allowed to overlap (this condition forces the
cosine of the minimal angle between the sphere centers to be at
most $1/3$, when the angle is centered at the midpoint between the
fixed balls). The largest such arrangement most likely consists of
$96$ balls.

To construct the code, consider three orthogonal tetrahedra in
$\R^9$.  Call the points in the first $v_1$, $v_2$, $v_3$, $v_4$,
in the second $v_5$, $v_6$, $v_7$, $v_8$, and in the third $v_9$,
$v_{10}$, $v_{11}$, $v_{12}$. Within each of these tetrahedra, all
inner products between distinct points are $-1/3$, and between
tetrahedra they are all $0$.  Call these tetrahedra the basic
tetrahedra.

The points $\pm v_1, \dots, \pm v_{12}$ will all be in the code,
and we will identify $72$ more points in it.  Each of the
additional points will have inner product $\pm 1/3$ with each of
$v_1,\dots,v_{12}$, and we will determine them via those inner
products.  Because $v_1+\dots+v_4 = v_5+\dots+v_8 =
v_9+\dots+v_{12} = 0$, the inner products with the elements of
each basic tetrahedron must sum to zero.  In particular, two must
be $1/3$ and the other two $-1/3$.  That restricts us to
$\binom{4}{2}=6$ patterns of inner products with each basic
tetrahedron, so there are $6^3=216$ points satisfying all the
constraints so far.  We must cut that number down by a factor of
$3$.

The final constraint comes from considering the inner products
between the new points.  A simple calculation shows that one can
reconstruct a point $x$ from its inner products with
$v_1,\dots,v_{12}$ via
$$
x = \frac{3}{4}\sum_{i=1}^{12} \langle x,v_i\rangle v_i,
$$
and inner products are computed via
$$
\langle x,y \rangle = \frac{3}{4} \sum_{i=1}^{12} \langle x, v_i
\rangle \langle y, v_i \rangle.
$$
In other words, if $x$ and $y$ have identical inner products with
one of the basic tetrahedra, that contributes $1/3$ to their own
inner product.  If they have opposite inner products with one of
the basic tetrahedra, that contributes $-1/3$.  Otherwise the
contribution is $0$.

The situation we wish to avoid is when $x$ and $y$ have identical
inner products with two basic tetrahedra, or opposite inner
products with both, and neither identical nor opposite inner
products with the third.  In that case, $\langle x,y \rangle = \pm
2/3$.

To rule out this situation, we assign elements of $\Z/3\Z$ to
quadruples by
$$
\pm (1/3,1/3,-1/3,-1/3) \mapsto 0,
$$
$$
\pm (1/3,-1/3,1/3,-1/3) \mapsto 1,
$$
and
$$
\pm (1/3,-1/3,-1/3,1/3) \mapsto -1.
$$
Consider the $72$ points with inner products $\pm 1/3$ with each
of $v_1,\dots,v_{12}$ such that exactly two inner products with
each basic tetrahedron are $1/3$ and furthermore the elements of
$\Z/3\Z$ coming from the inner products with the basic tetrahedra
sum to $0$.  Given any two such points, if they have identical or
opposite inner products with two basic tetrahedra, then the same
must be true with the third. Thus, we have constructed $24+72=96$
points in $\R^9$ such that all the inner products between them are
$\pm 1$, $\pm 1/3$ or $0$.

\begin{table}
\begin{center}
\begin{tabular}{ccc}
Vertices & Automorphisms & Orbit size \\
\hline $9$ & $16$ & $27648$\\
$9$ & $48$ & $13824$\\
$9$ & $48$ & $4608$\\
$9$ & $96$ & $18432$\\
$9$ & $1440$ & $4608$\\
$12$ & $1024$ & $864$\\
$12$ & $31104$ & $512$\\
$16$ & $10321920$ & $18$\\
\\ % Bad to skip line?
\end{tabular}
\end{center}
\caption{Facets of the convex hull of the configuration of $96$
points in $\R^9$, modulo the action of the symmetry
group.}\label{table:9-96facets}
\end{table}

The facets of this code form eight orbits under the action of its
symmetry group; they are listed in Table~\ref{table:9-96facets}. The
most interesting facets are those with $16$ vertices, which form
regular cross polytopes.  These facets and the two orbits with $12$
vertices all correspond to deep holes.

\subsection{Distribution of energy levels}

Typically, there are many local minima for harmonic energy. One
intriguing question is how the energies of the local minima are
distributed.  For example, Table~\ref{table:energy} shows the
thirty lowest energies obtained in $2\cdot10^5$ trials with $120$
points in $\R^4$, together with how often they occurred.  The
regular $600$-cell is the unique universal optimum (with energy
$5395$), but we found $5223$ different energy levels.  This table
is probably not a complete list of the lowest thirty energies
(five of them occurred only once, so it is likely there are more
to be found), but we suspect that we have found the true lowest
ten.

\begin{table}
\begin{center}
\begin{tabular}{cc|cc}
Energy & Frequency & Energy & Frequency\\
\hline $5395.000000$ & $186418$ & $5402.116636$ & $1$\\
$5398.650556$ & $4393$ & $5402.152619$ & $1$\\
$5398.687876$ & $2356$ & $5402.213231$ & $2$\\
$5400.842726$ & $18$ & $5402.366164$ & $1$\\
$5400.880057$ & $149$ & $5402.922701$ & $1$\\
$5400.890460$ & $47$ & $5403.091064$ & $111$\\
$5400.894513$ & $26$ & $5403.115123$ & $1$\\
$5400.928674$ & $25$ & $5403.129076$ & $108$\\
$5400.936106$ & $41$ & $5403.271100$ & $66$\\
$5400.940237$ & $28$ & $5403.319898$ & $157$\\
$5400.940550$ & $7$ & $5403.326719$ & $84$\\
$5400.943094$ & $38$ & $5403.347209$ & $24$\\
$5402.029556$ & $7$ & $5403.455701$ & $7$\\
$5402.088248$ & $3$ & $5403.462898$ & $8$\\
$5402.093726$ & $10$ & $5403.488923$ & $4$\\
\\
\end{tabular}
\end{center}
\caption{Thirty lowest harmonic energies observed for local minima
with $120$ points on $S^3$ ($2 \cdot 10^5$ trials, $5223$
different energy levels observed).} \label{table:energy}
\end{table}

The most remarkable aspect of Table~\ref{table:energy} is the
three gaps in it.  There are huge gaps from $5395$ to $5398.65$,
from $5398.69$ to $5400.84$, and from $5400.95$ to $5402.02$. Each
gap is far larger than the typical spacing between energy levels.
Perhaps one of these gaps contains some rare local minima, but
they appear to be real gaps.

\begin{figure}
\begin{center}
\includegraphics{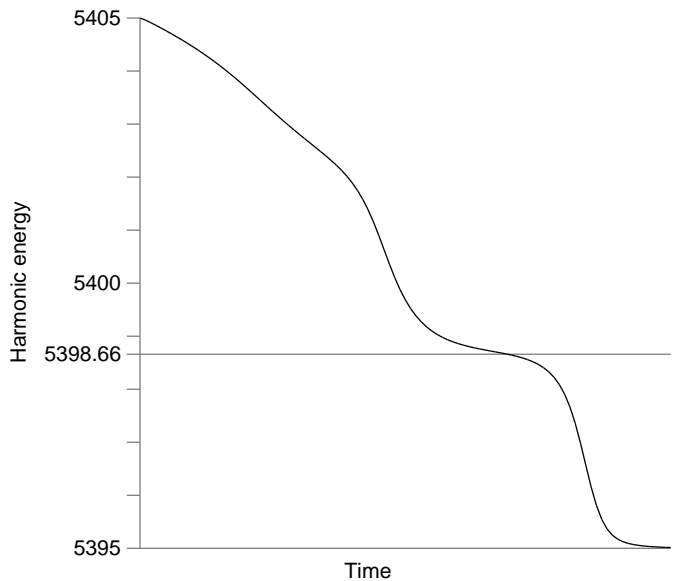}
\end{center}
\caption{Energy as a function of time under gradient descent.}
\label{fig:progress}
\end{figure}

What could cause such gaps?  We do not have a complete theory, but
we believe the gaps reflect bottlenecks in the process of
constructing the code by gradient descent.
Figure~\ref{fig:progress} is a graph of energy as a function of
time for gradient descent, starting at a random configuration of
$120$ points in $S^3$.  (The figure represents a single run of the
optimization procedure, so it should be viewed as a case study,
not a statistical argument.) The graph begins when the energy has
just reached $5405$ and ends once convergence to the universal
optimum is apparent. Of course the convergence is monotonic, but
its speed varies dramatically. The rate of decrease is slowest
near energy $5398.66$, which is indicated by a horizontal line. We
do not believe it could be a coincidence that that is very nearly
the energy of the two lowest-energy local minima in
Table~\ref{table:energy}: the slowdown probably occurs because of
a bottleneck.  More precisely, in order to achieve the ground
state the system must develop considerable large-scale order and
symmetry.  Probably short-range order develops first and then
slowly extends to long-range order. During this process there may
be bottlenecks in which different parts of the system must come
into alignment with each other.  The local minima correspond to
the rare cases in which the system gets stuck in the middle of a
bottleneck, but even when it does not get stuck it still slows
down.

It is not completely clear why the gaps are separated by several
energy levels that are surprisingly close to each other.  The most
likely explanation is that there are several slightly different
ways to get stuck during essentially the same bottleneck, but we
have no conceptual understanding of what these ways are. It would
be very interesting to have a detailed theory of this sort of
symmetry breaking.

\section{Conjectured universal optima}
\label{section:newunivopt}

\subsection{$40$ points in $\R^{10}$}

The $40$-point code from Table~\ref{table:conjectures} consists of
$40$ points on the unit sphere in $\R^{10}$.  The only inner
products that occur are $1$, $-1/2$, $-1/3$, $0$, and $1/6$; each
point has these inner products with $1$, $8$, $3$, $4$, and $24$
points, respectively. Grouping pairs of points according to their
inner product yields a $4$-class association scheme, which Bannai,
Bannai, and Bannai \cite{BBB} have recently shown is uniquely
determined by its intersection numbers.

The $40$ points form $10$ regular tetrahedra, and more
specifically $5$ orthogonal pairs of tetrahedra.  That accounts
for all the inner products of $1$, $-1/3$, and $0$.  Each point
has inner product $-1/2$ with one point in each of the other $9$
tetrahedra, except for the tetrahedron orthogonal to the one
containing it. All remaining inner products are $1/6$. The
configuration is chiral (i.e., not equivalent to any reflection of
itself under the action of $SO(10)$) and has a symmetry group of
order $1920$. Specifically, the symmetry group is the semidirect
product of the symmetric group $S_5$ with the subgroup of
$(\Z/2\Z)^5$ consisting of all vectors that sum to zero, where
$S_5$ acts by permuting the five coordinates.

We conjecture that this code is the unique $40$-point code in
$\R^{10}$ with maximal inner product $1/6$, but that appears
difficult to prove.  For example, it is an even stronger assertion
than optimality as a spherical code.  We are unaware of any
occurrence of this code in the published literature, but it
appears in Sloane's online tables \cite{S} with the annotation
that it was found by Smith and ``beautified'' by Conway and Sloane
(in the sense of recognizing it using elegant coordinates).  It
also appears in Hovinga's online report \cite{H}.

Conway, Sloane, and Smith construct the code as an explicit list
of $40$ vectors with entries $-1/\sqrt{6}$, $0$, or $1/\sqrt{6}$.
Here we explain the combinatorics that underlies this
construction. That may well be how Conway and Sloane beautified
the code, but \cite{S} presents no details about the construction
beyond the list of vectors.

Consider the $10$-dimensional vector space $V$ spanned by the
orthonormal basis vectors $v_{\{i,j\}}$ for all two-element
subsets $\{i,j\} \subset \Z/5\Z$ with $i \ne j$. Say a vector in
$V$ has type $i$ if for every basis vector $v_S$, its coefficient
vanishes if and only if $i \in S$.  Each vector in the code will
have type $i$ for some $i \in \Z/5\Z$, and the six nonzero
coefficients will equal $\pm 1/\sqrt{6}$.  Given such a vector,
define a graph on the vertex set $(\Z/5\Z) \setminus \{i\}$ by
connecting $j$ to $k$ if the coefficient of $v_{\{j,k\}}$ is
$-1/\sqrt{6}$.  If the graph has $e$ edges and vertex $j$ has
degree $d_j$, then the vector is in the code if and only if
$d_{i+1} \equiv d_{i+2} \equiv e \pmod{2}$ and $d_{i-1} \equiv
d_{i-2} \not\equiv e \pmod{2}$ (where $i$ is the type of the
vector).

The facet structure of this code seems surprisingly unilluminating
compared to the others we have analyzed.  There are $24$ orbits of
facets: $21$ orbits of irregular simplices (with symmetry groups
ranging from $2$ to $768$ in size), one orbit of $11$-vertex
facets, and two orbits of $12$-vertex facets.  The most
symmetrical facets are those in one of the two orbits of
$12$-vertex facets. They have $3072$ symmetries and are given by
the orthogonal union of three identical irregular tetrahedra in
$\R^3$ (with $8$ symmetries).

Our confidence in this code's universal optimality is based on detailed
numerical experiments. One reason a configuration could be universally
optimal is that it has no competitors (i.e., except in degenerate cases
there are no other local minima). That is not true for the $40$-point
code, but there are remarkably few competitors.  In particular, it
appears to have only one ``serious'' competitor. We have found five
other families of local minima, but four of them are rare and never
seem to come close to beating our conjectured universal optimum.  The
best experimental evidence we can imagine for universal optimality
would be to describe explicitly each competing family that has been
observed and prove that it never contains the global energy minimum.
That might be possible for this code, but we have not completed it. The
four rare families are sufficiently complicated that we have not
analyzed them explicitly (under the circumstances it did not seem worth
the effort). However, we have a complete description of the serious
competitor.

That family depends on a parameter $\alpha$ that must satisfy $0 <
\alpha^2 \le 1/27$.  The configuration always contains a fixed
$16$-point subset with the following structure.  If we call the
$16$ points $w_{i,j}$ with $i,j \in \{1,2,3,4\}$, then
$$
\langle w_{i,j}, w_{k,\ell}\rangle = \begin{cases} 1 & \textup{if
$(i,j) = (k,\ell)$},\\
-1/3 & \textup{if $i=k$ or $j =\ell$ but not both, and}\\
1/9 & \textup{otherwise.}
\end{cases}
$$
In other words, if we arrange the points in a $4 \times 4$ grid,
then the rows and columns are regular tetrahedra and all other
inner products are $1/9$.  To construct such a configuration, take
the tensor product in $\R^3 \otimes_\R \R^3 = \R^9$ of two regular
tetrahedra in $\R^3$.

To describe the remaining $24$ points, we will specify their inner
products with the first $16$ points.  That will determine their
projections into the $9$-dimensional subspace containing the $16$
points, so the only additional information needed to pin them down
will be whether they are above or below that hyperplane (relative
to some orientation).

Each of the $24$ points will have inner product $-3\alpha$ with
four of the first $16$ points and $\alpha$ with each of the
others.  The only constraint is that it must have inner product
$-3\alpha$ with exactly four points, one in each of the eight
tetrahedra (i.e., one in each row and column of the $4\times 4$
grid).  The $4\times 4$ grid exhibits a one-to-one correspondence
with permutations of four elements, so there are $4!=24$ ways to
satisfy these constraints.  The points corresponding to even
permutations will be placed above the $9$-dimensional hyperplane,
and those corresponding to odd permutations will be placed below
it. This construction yields a $40$-point code in $\R^{10}$
provided that $0 < \alpha^2 \le 1/27$.  (When $\alpha=0$ some
points coincide, and when $\alpha^2 > 1/27$ the inner products
cannot be achieved by unit vectors.)

We have not proved that the codes in this family never improve on
the conjectured universal optimum, but we are confident that it is
true. The best spherical code in the family occurs when $\alpha =
(\sqrt{109}-1)/54 = 0.1748\dots$; in this special case, $\alpha$
is also the maximal inner product, which is quite a bit larger
than the maximal inner product $1/6$ in the conjectured universal
optimum. That implies that when $k$ is sufficiently large, the
conjectured optimum is better for the potential function $f(r) =
(4-r)^k$ (because the energy is dominated asymptotically by the
contribution from the minimal distance).  In principle one could
verify the finitely many remaining values of $k$ by a finite
computation.  We have done enough exploration to convince
ourselves that it is true, but we have not found a rigorous proof.

\subsection{$64$ points in $\R^{14}$}
\label{subsection:64}

The simplest construction we are aware of for the $64$-point
configuration in $\R^{14}$ uses the Nordstrom-Robinson binary code
\cite{NR,G}.  Shortening that code twice yields a binary code of
length $14$, size $64$, and minimal distance $6$, which is known
to be unique (see \cite[pp.~74--75]{MS}). One can view it as a
subset of the cube $\{-1,1\}^{14}$ instead of $\{0,1\}^{14}$. Then
after rescaling by a factor of $1/\sqrt{14}$ to yield unit
vectors, this code is the $64$-point configuration in $\R^{14}$
that we conjecture is universally optimal. The same process with
less shortening yields the codes from Tables~\ref{table:balanced}
and~\ref{table:unresolved} with $128$ points in $\R^{15}$ and
$256$ points in $\R^{16}$;
% fix me
the $64$-point and $128$-point codes have previously appeared in
Appendix~D of \cite{EZ} via the same approach, as conjectures for
optimal spherical codes.

This construction makes some of the facet structure of the code
clear.  There are $28$ facets with $32$ vertices that come from
the facets of the cube containing the code.  In fact, the code is
obtained from exactly the same antiprism construction as described
in Subsection~\ref{subsec:facet} (the vertices of these facets
point towards deep holes of the opposite facets).  There are also
$66$ other orbits of facets under the action of the symmetry
group, but those orbits seem to be less interesting.

An alternative construction of the code is by describing its Gram
matrix explicitly.  As mentioned above, this construction amounts
to forming an association scheme by taking $t=1$ in Theorem~2 and
Proposition~7(i) in \cite{dCvD}, and then performing a spectral
embedding.  Bannai, Bannai, and Bannai \cite{BBB} have shown that
this association scheme is uniquely determined by its intersection
numbers.

More concretely, the points correspond to elements of $\F_8^2$,
where $\F_8$ is the finite field of order $8$, and the inner
products are determined by Table~\ref{table:gram64}.  The Gram
matrix has $14$ eigenvalues equal to $32/7$ and the others equal
to $0$, so it is indeed the Gram matrix of a $14$-dimensional
configuration.

\begin{table}
\begin{tabular}{cc}
Inner product & Condition\\
\hline
$1$ & $(x_1,x_2)=(y_1,y_2)$\\
$-1/7$ & $x_1=y_1$ but $x_2 \ne y_2$\\
$-3/7$ & $x_1 \ne y_1$ and $x_2+y_2 \in \{ (x_1+y_1)^3, x_1y_1(x_1+y_1)\}$\\
$1/7$ & otherwise\\
\\
\end{tabular}
\caption{The inner products for the $64$-point code.}
\label{table:gram64}
\end{table}

Unfortunately, this $64$-point code has many competitors.  We have
found over two hundred local minima for harmonic energy and expect
that the total number is much larger.  That makes it difficult to
imagine an ironclad experimental argument for universal
optimality.  We suspect that the code is universally optimal for
two reasons: we have failed to find any counterexample, and during
the process no competitor came close enough to worry us.  (By
contrast, in many cases in which one can disprove universal
optimality, one finds worrisomely close competitors before
tweaking the construction to complete the disproof.)  However, we
realize that the evidence is far from conclusive.

\section{Balanced, irreducible harmonic optima}
\label{section:balanced}

In this section we briefly describe each of the configurations in
Tables~\ref{table:balanced} and~\ref{table:unresolved}.

% The \medskip commands that follow are a kludge:

\medskip

\begin{description}

\item[$32$ points in $\R^3$]

The union of a regular icosahedron and its dual dodecahedron.

\medskip

\item[$10$ points in $\R^4$]

Both configurations are described in
Subsection~\ref{subsec:diplo}.

\medskip

\item[$13$ points in $\R^4$]

In $\C^2$ with the inner product $\langle x,y \rangle =
\mathop{\textup{Re}} \bar x ^t y$, the points are
$(\zeta/\sqrt{2},\zeta^5/\sqrt{2})$,  where $\zeta$ runs over all
$13$-th roots of unity.  This code was discovered by Sloane,
Hardin, and Cara \cite{SHC}.

For a less compact description, view $\R^4$ as the orthogonal
direct sum of two planes, and let $R$ be the operation of rotating
the first plane by $2\pi/13$ and the second by five times that
angle. The unit sphere in $\R^4$ contains the direct product of
the circles of radius $1/\sqrt{2}$ in the two planes, and the
$13$-point code is the orbit of a point in this direct product
under the group generated by $R$.

The factor of $5$ is special because $5^2 \equiv -1 \pmod{13}$. In
particular, $R^8$ rotates the second plane by $2\pi/13$ and the
first by $8$ times that amount, which is the same as $5$ times it
in the opposite direction.  In other words, the two planes play
the same role, if one ignores their orientations.  Only the square
roots of $-1$ modulo $13$ (or, trivially, the square roots of $1$)
have that property.

\begin{figure}
\begin{center}
\includegraphics{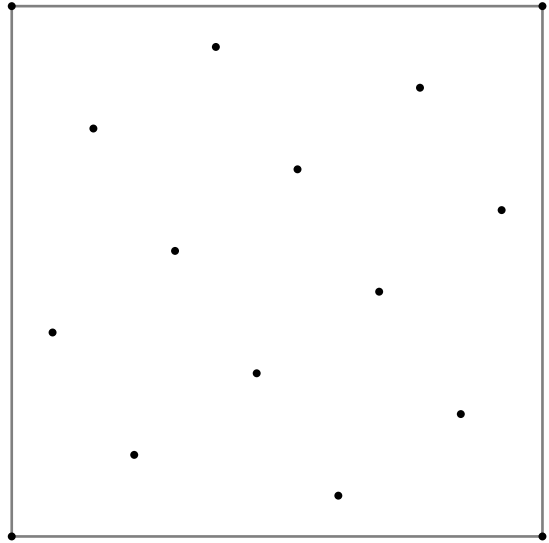}
\end{center}
\caption{The $13$-point harmonic optimum in $S^3$, drawn on a
torus by plotting $(\phi,\psi)$ for the point
$(e^{i\phi}/\sqrt{2},e^{i\psi}/\sqrt{2})$.} \label{fig:torus13}
\end{figure}

For all the points in this configuration, both complex coordinates
have absolute value $1/\sqrt{2}$.  Therefore the configuration is
contained in a flat two-dimensional torus sitting inside $S^3$
(namely, the product of the circles of radius $1/\sqrt{2}$ in the
two complex coordinate axes).  Figure~\ref{fig:torus13} shows the
complex phases of the $13$ points.  For each $N$ we can ask whether
there is an $N$-point harmonic optimum in $S^3$ that is contained
in such a torus.  For $1 \le N \le 10$ it seems that there is.  We
conjecture that $N=13$ is the only larger value of $N$ for which
this happens.

For other potential functions, similar phenomena can occur in more
cases. For example, for the logarithmic potential function $f(r) =
-\log r$, Jaron Lanier has conjectured in a private communication
that the minimal-energy configuration of $11$ points on $S^3$ also
lies on a flat torus.  In $\C^2$, the points are
$(\alpha\zeta,\sqrt{1-\alpha^2}\zeta^4)$, where $\zeta$ runs over
all $11$-th roots of unity and $\alpha$ is the unique root between
$0$ and $1$ of
$$
5 \alpha^8 - 36 \alpha^6 + 51 \alpha^4 - 4 \alpha^2 - 7= 0.
$$

\medskip

\item[$15$ points in $\R^4$]

Described in Subsection~\ref{subsec:gramanalysis}.

\medskip

\item[$24$ points in $\R^4$]

The regular $24$-cell (equivalently, the $D_4$ root system).

\medskip

\item[$48$ points in $\R^4$]

Described in Subsection~\ref{subsec:48in4}.

\medskip

\item[$21$ points in $\R^5$]

The edge midpoints and face centers of a regular simplex (rescaled
to lie on the same sphere).

\medskip

\item[$32$ points in $\R^5$]

Start with a regular simplex $v_1,\dots,v_6$ in $\R^5$.  The $32$
points will be these six points and their negatives, along with
$20$ points determined  as follows by their inner products with
$v_1,\dots,v_6$.  Each will have inner product $\pm 1/\sqrt{5}$
with each of $v_1,\dots,v_6$, with three plus signs and three
minus signs.  There are $\binom{6}{3}=20$ ways to choose these
signs.

\medskip

\item[$2n+2$ points in $\R^n$ (for $n \ge 6$)]

Described in Subsection~\ref{subsec:diplo}.

\medskip

\item[$42$ points in $\R^6$]

The edge midpoints of a regular simplex and their antipodes.

\medskip

\item[$44$ points in $\R^6$]

This code contains plus or minus an orthonormal basis of $\R^6$
together with the $32$ vectors whose coordinates with respect to
that basis are $\pm 1/\sqrt{6}$ and where an even number of minus
signs occur.  It other words, it consists of a cross polytope and
a hemicube within the cube dual to the cross polytope.  This code
was previously conjectured to be an optimal spherical code (see
Table~D.6 in \cite{EZ}).

\medskip

\item[$126$ points in $\R^6$]

The union of the minimal vectors of the $E_6$ and $E_6^*$ lattices
(rescaled to lie on the same sphere).  Equivalently, one can
project the $E_7$ root system orthogonally to a minimal vector in
$E_7^*$, followed by rescaling as in the first construction.

\medskip

\item[$78$ points in $\R^7$]

Like the $44$ points in $\R^6$, this code consists of a cross
polytope and a hemicube within the cube dual to the cross
polytope.

\medskip

\item[$148$ points in $\R^7$]

The points are all the permutations of
$$
\frac{(\pm1,\pm1,0,0,0,0,0)}{\sqrt{2}}
$$
and the hemicube consisting of the points
$$
\frac{(\pm1,\pm1,\pm1,\pm1,\pm1,\pm1,\pm1)}{\sqrt{7}}
$$
that have an even number of minus signs.  This is a
seven-dimensional analogue of a construction of the $E_8$ root
system, but it is less symmetric, because the two types of points
form distinct orbits under the action of the symmetry group.  This
code was previously conjectured to be an optimal spherical code
(see Table~D.7 in \cite{EZ}).

\medskip

\item[$182$ points in $\R^7$]

The union of the minimal vectors of the $E_7$ and $E_7^*$ lattices
(rescaled to lie on the same sphere).  Equivalently, one can
project the $E_8$ root system orthogonally to any root; $126$
roots are unchanged, $112$ project to $56$ nonzero points, and $2$
project to the origin.  Rescaling the nonzero projections to lie
on the unit sphere yields the $182$-point configuration.

\medskip

\item[$72$ points in $\R^8$]

The edge midpoints of a regular simplex and their antipodes.  This
code was previously conjectured to be an optimal spherical code
(see Table~D.8 in \cite{EZ}).

\medskip

\item[$96$ points in $\R^9$]

Described in Subsection~\ref{subsec:96in9}.

\medskip

\item[$42$ points in $\R^{14}$]

This code consists of seven disjoint, five-dimensional regular
simplices, whose vertices have inner products $-1/5$ with each
other. Each point has inner product $-1/2$ with a unique point in
each simplex other than the one containing it, and all other inner
products between points in different simplices are $1/10$. Grouping
pairs of points according to their inner product yields a
three-class association scheme that can be derived from the
Hoffman-Singleton graph (see Subsection~5.1 of \cite{vD}). Let $G$
be the Hoffman-Singleton graph, and let $H$ be its second
subconstituent. In other words, take any vertex $v$ in $G$, and let
$H$ be the vertices not equal to or adjacent to $v$. The vertices
in $H$ correspond to points in the $14$-dimensional configuration,
with inner product $-1/2$ between adjacent vertices, $-1/5$ between
non-adjacent vertices with no common neighbor, and $1/10$ between
non-adjacent vertices with one common neighbor.  This code was
previously conjectured to be an optimal spherical code (see
Table~D.14 in \cite{EZ}).

\medskip

\item[$128$ points in $\R^{15}$]

Described in Subsection~\ref{subsection:64}.

\medskip

\item[$256$ points in $\R^{16}$]

Described in Subsection~\ref{subsection:64}.

\end{description}

\section{Challenges}

We conclude with a list of computational and theoretical
challenges:

\begin{enumerate}
\item How often do $196560$ randomly chosen points on $S^{23}$
converge to the Leech lattice minimal vectors under gradient
descent for harmonic energy?  For $240$ points on $S^7$, one
frequently obtains the $E_8$ root system ($855$ times out of
$1000$ trials), and the Higman-Sims configuration of $100$ points
on $S^{21}$ occurs fairly often ($257$ times out of $1000$
trials); by contrast, the universal optimum with $112$ points on
$S^{20}$ occurs rarely (once in $1000$ trials). This approach
could be an intriguing construction of the Leech lattice, but we
have no intuition for how likely it is to work.

\item What are the potential energy barriers that separate the local
minima in Table~\ref{table:27} (or any other case)?  In other
words, if one continuously transforms one configuration into
another, how low can one make the greatest energy along the path
connecting them? The lowest possible point of greatest energy will
always be a saddle point for the energy function.

\item How many harmonic local minima are there for $120$ points on $S^3$
(Table~\ref{table:energy}), or even $64$ points on $S^{13}$? Is
the number small enough that one could conceivably compile a
complete list?  Is the list for $27$ points on $S^5$ in
Table~\ref{table:27} complete?

\item For large numbers of points, what can one say (experimentally,
heuristically, or rigorously) about the distribution of energy
levels for local minima, or about the gaps in the distribution?
\end{enumerate}

\section*{Acknowledgements}

We thank Eiichi Bannai, Christian Borgs, John Conway, Edwin van
Dam, Charles Doran, Noam Elkies, Florian Gaisendrees, Robert
Griess, Abhinav Kumar, Jaron Lanier, James Morrow, Frank
Stillinger, and Salvatore Torquato for helpful discussions.

\appendix
\section{Local non-optimality of diplo-simplices}
\label{appendix:diplo}

In this appendix we prove that (above dimension two)
diplo-simplices are not even locally optimal as spherical codes.
The motivation behind the calculations is the case of $n=3$, where
the diplo-simplex is a cube.  One can improve this spherical code
by rotating one face relative to the opposite face while moving
them slightly closer together.  To carry out this approach in
higher dimensions, we must understand the faces of the
diplo-simplex (which is in general quite different from the
hypercube).   The diplo-simplex in $\R^n$ is the orthogonal
projection of the vertices of the cross polytope in $\R^{n+1}$
onto the hyperplane on which the sum of the coordinates is zero.
Its dual polytope is therefore the cross section of the hypercube
in $\R^{n+1}$ by that hyperplane. The vertices of the cross
section differ depending on whether $n$ is odd or even:
$(1,\dots,1,-1,\dots,-1)$ is a vertex when $n+1$ is even, and
$(1,\dots,1,0,-1,\dots,-1)$ is when $n+1$ is odd.  That leads to
the following description of the diplo-simplex, based on
identifying its faces using this correspondence.

Let $v_1,\dots,v_k$ and $w_1,\dots,w_k$ be unit vectors forming
two orthogonal, $(k-1)$-dimensional regular simplices, and let $t$
be a unit vector orthogonal to both simplices.  The
$(2k-1)$-dimensional diplo-simplex consists of the points
$$
\pm(\alpha v_i + \sqrt{1-\alpha^2}t)
$$
and
$$
\pm(\alpha w_i + \sqrt{1-\alpha^2}t),
$$ where $\alpha = \sqrt{(2k-2)/(2k-1)}$. To
improve its minimal angle, use the points
$$
\alpha v_i + \sqrt{1-\alpha^2}t,
$$
$$
\alpha w_i + \sqrt{1-\alpha^2}t,
$$
$$
-\big(\alpha (\beta v_i + \sqrt{1-\beta^2} w_i) +
\sqrt{1-\alpha^2}t\big),
$$
and
$$
-\big(\alpha (\sqrt{1-\beta^2} v_i - \beta w_i) +
\sqrt{1-\alpha^2}t\big),
$$
with $\beta$ a small positive number and $\alpha$ slightly greater
than $\sqrt{(2k-2)/(2k-1)}$.  When $k=1$ choosing $\alpha$ and
$\beta$ optimally leads to the optimal $8$-point code in $\R^3$,
but when $k=3$ Sloane's tables \cite{S} show that this approach is
suboptimal.

The situation is slightly more complicated for even-dimensional
diplo-simplices. Again, let $v_1, \dots, v_k$ and $w, \dots, w_k$ be
unit vectors forming two orthogonal $(k-1)$-dimensional simplices. Let
$t$ and $z$ be two unit vectors orthogonal to both simplices and each
other. The $2k$-dimensional diplo-simplex consists of the points $z$,
$-z$,
$$
\pm (\alpha v_i + \beta z + \sqrt{1-\alpha^2-\beta^2}t),
$$
and
$$
\pm (\alpha w_i - \beta z +
\sqrt{1-\alpha^2-\beta^2}t),
$$ where
$\alpha=\sqrt{(2k+1)(2k-2)}/(2k)$ and $\beta=1/(2k)$. Note that
$\alpha$ vanishes when $k=1$ and our construction below will not
work. Indeed the hexagon, which is the diplo-simplex in the plane,
is an optimal spherical code. To improve the minimal angle for $k
\geq 2$ use the points $z$, $-z$,
$$
\alpha v_i + \beta z + \sqrt{1-\alpha^2-\beta^2}t,
$$
$$
\alpha w_i - \beta z + \sqrt{1-\alpha^2-\beta^2}t,
$$
$$
-\big(\alpha(\gamma v_i + \sqrt{1-\gamma^2}w_i)-\beta z +
\sqrt{1-\alpha^2 - \beta^2}t\big),
$$
and
$$
-\big(\alpha(\sqrt{1-\gamma^2}v_i - \gamma w_i)+\beta z +
\sqrt{1-\alpha^2 - \beta^2}t\big),
$$
with $\gamma$ a small positive number, $\alpha$ slightly larger
than $\sqrt{(2k+1)(2k-2)}/(2k)$ and $\beta$ slightly less than
$1/(2k)$. The numbers $\alpha$ and $\beta$ should be chosen such
that $\alpha^2 + 2 \beta^2$ increases from its original value of
$(2k-1)/(2k)$.

\end{document}